\documentclass[a4paper,fleqn,11pt]{elsarticle}

\usepackage{amssymb,amsthm,amsmath}    
\usepackage{a4wide}
\usepackage[english]{babel}
\usepackage{csquotes}
\usepackage{thmtools} %

\usepackage{graphicx} 
\usepackage{tikz,pgfplots} 
    \pgfplotsset{compat=1.14}
    \usepgfplotslibrary{fillbetween}
    \pgfdeclarelayer{ft}
    \pgfdeclarelayer{bg}
    \pgfsetlayers{bg,main,ft}
\usepackage{hyperref} 
\hypersetup{colorlinks,linkcolor={red!50!black},citecolor={blue!50!black},urlcolor={blue!80!black}}
    \addto\extrasenglish{%
    }
\usepackage{enumitem}
\usepackage{subcaption}
\usepackage{natbib}
\usepackage{doi}

\usepackage{accents}
\usepackage{etoolbox}
\usepackage{array}
\newcolumntype{C}[1]{>{\centering\arraybackslash $}p{#1}<{$}}

\theoremstyle{plain}
\newtheorem{theorem}{Theorem}[section]
\newtheorem{lemma}[theorem]{Lemma}
\newtheorem{proposition}[theorem]{Proposition}

\newtheorem{problem}[theorem]{Problem}

\theoremstyle{definition}
\newtheorem{definition}[theorem]{Definition}
\newtheorem{remark}[theorem]{Remark}
\newtheorem{example}[theorem]{Example}

\newcommand{\Z}{\mathbb{Z}} 
\newcommand{\N}{\mathbb{N}} 
\newcommand{\R}{\mathbb{R}} 

\newcommand{\alphabet}{\Sigma}
\newcommand{\words}{\alphabet^*} 
\newcommand{\eps}{\varepsilon} 

\DeclareMathOperator{\pref}{pref} 
\newcommand{\freq}[2][{}]{\operatorname{freq}_{#1}(#2)} 
\newcommand{\supfreq}[2][{}]{\overline{\operatorname{freq}}_{#1}(#2)} 
\newcommand{\inffreq}[2][{}]{\underline{\operatorname{freq}}_{#1}(#2)} 

\newcommand{\infw}[1]{\mathbf{#1}}

\newcommand{\lang}[2][{}]{\mathcal{L}_{#1}(#2)}

\newcommand{\soc}[1]{\Omega(#1)}
\newcommand{\abclsr}[1]{\mathcal{A}(#1)}

\newcommand{\complfunction}[1]{\mathcal{P}_{#1}}

\newcommand{\gap}{\hskip2pt}

\begin{document}

    \title{On Abelian Closures of Infinite Non-binary Words}
    \author[1]{Juhani Karhum\"aki}
		\ead{karhumak@utu.fi}
		
		\author[2,3]{Svetlana Puzynina
		\corref{cor1}}
		\ead{s.puzynina@gmail.com}
		
		\author[4]{Markus A. Whiteland\corref{cor1}
		}
		\ead{mawhit@mpi-sws.org}
		
	\address[1]{Department of Mathematics and Statistics, FI-20014 University of Turku, Finland}
		\address[2]{St. Petersburg State University, Russia}
		\address[3]{Sobolev Institute of Mathematics, Russia}
		\address[4]{Max Planck Institute for Software Systems, Saarland Informatics Campus,
		Saarbr{\"u}cken, Germany}

\begin{abstract}
Two finite words $u$ and $v$ are called abelian equivalent
if each letter occurs equally many times in both $u$ and $v$. The
abelian closure $\abclsr{\infw{x}}$ of 
an infinite word
$\infw{x}$ is the set of infinite words $\infw{y}$ such that, for each
factor $u$ of $\infw{y}$, there exists a factor $v$ of $\infw{x}$
which is abelian
equivalent to $u$. 
The notion of an abelian closure gives a characterization of
Sturmian words: among uniformly recurrent binary words, periodic and aperiodic Sturmian
words are exactly those words for which $\abclsr{\infw{x}}$
equals the shift orbit closure $\soc{\infw{x}}$. Furthermore,
for an aperiodic binary word that is not Sturmian, its abelian closure
contains infinitely many minimial subshifts.
In this paper we consider the abelian closures of well-known families of
non-binary words, such as balanced words and minimal complexity words.
We also consider abelian closures of general subshifts and make some initial
observations of their abelian closures and pose some related open questions.
\end{abstract}

\maketitle


\section{Introduction}

Let ${\infw x}\in \Sigma^{\mathbb{N}}$ be an infinite word over an
alphabet $\Sigma$. We define the \emph{language of $\infw x$},
denoted by $\mathcal L({\infw x})$, as the set of factors of ${\infw
x}$, i.e., blocks of consecutive letters of ${\infw x}$. A
\emph{subshift} $\soc{\infw x}$ generated by an infinite word
${\infw x}$ can be defined as the set of infinite words whose
languages are included in $\mathcal L({\infw x})$: $\soc{\infw x}
=\{{\infw y}\in \Sigma^{\mathbb{N}}\colon \mathcal L({\infw
y})\subseteq \mathcal L({\infw x}) \}.$ In this paper, we consider
an abelian version of the notion of a subshift. Two finite words
$u$ and $v$ are called \emph{abelian equivalent}, denoted by
$u\sim_{\text{ab}} v$, if each letter occurs equally many times in
both $u$ and $v$. Various abelian properties of words have been
actively studied recently, e.g., abelian complexity, abelian
powers, abelian periods, etc.
\cite{DBLP:journals/eatcs/ConstantinescuI06,PUZYNINA2013390,DBLP:journals/jlms/RichommeSZ11}.
We define the \emph{abelian closure} $\abclsr{\infw x}$ of an
infinite word $\infw x$ as the set of infinite words $\infw y$ such
that, for each factor $u$ of $\infw y$, there exists a factor $v$
of $\infw x$ with $u\sim_{\text{ab}}v$. Clearly, $\soc{\infw x}
\subseteq \abclsr{\infw x}$ for any word $\infw x$.

We start with two examples showing completely different structure
of abelian closures: Sturmian words and the Thue--Morse word.
Sturmian words can be defined as infinite aperiodic words which
have $n+1$ distinct factors for each length $n$. They admit
various characterizations; in particular, they are exactly the
aperiodic binary balanced words (i.e., the numbers of occurrences
of $1$ in factors of the same length differ by at most $1$). It is
not hard to see that, for a Sturmian word $\infw x$, $\soc{\infw
x} = \abclsr{\infw x}$ (so, the abelian closure is small, contains
only its subshift). Indeed, due to balance there are exactly two
abelian classes of factors of each length. Therefore, any word
$\infw{y}\in \abclsr{\infw{x}}$ must be balanced. Further, the
frequencies of letters of ${\infw y}$ are uniquely defined by
$\abclsr{\infw x}$. Thus ${\infw y}$ is Sturmian with the same
letter frequencies as $\infw{x}$, i.e., $\infw{y} \in \soc{\infw{x}}$.
In fact, the property $\soc{\infw x} = \abclsr{\infw x}$ characterizes
Sturmian words among uniformly recurrent binary words
(see \autoref{th:St}).

The Thue--Morse word $\infw{TM}=011010011001\cdots$ can be defined
as the fixed point starting with $0$ of the morphism $\mu:0\mapsto
01$, $1\mapsto 10$. For odd lengths $\infw{TM}$ has two abelian
factors, and for even lengths three. Further, the number of
occurrences of $1$ in each factor differs by at most $1$ from half
of its length \cite{DBLP:journals/jlms/RichommeSZ11}. It is easy
to see that any factor of any word in $\{\eps,0,1\}\cdot
\{01,10\}^{\N}$ has the same property, i.e., $\{\eps,0,1\}\cdot
\{01,10\}^{\N}\subseteq \abclsr{\infw{TM}}$. In fact, equality
holds: $\abclsr{\infw{TM}}=\{\eps,0,1\}\cdot \{01,10\}^{\N}$
(so, the abelian closure of $\infw{TM}$ is huge compared to
its shift orbit closure). 
Indeed, let $\infw x\in \mathcal
A({\infw{TM}})$. Then $\infw x$ has blocks of each letter of length
at most 2 (since there are no factors $000$ and $111$). Moreover,
between two consecutive occurrences of $00$ there must occur $11$,
and vice versa (otherwise we have a factor $00(10)^n0$, where the
number of occurrences of $1$ differs by more than $1$ from half of
its length). Clearly, such a word is in $\{\eps,0,1\}\cdot
\{01,10\}^{\N}.$ In \cite{DBLP:journals/corr/abs-2008-08125}, we
show that this fact can be generalized to all  binary words: In
fact, each binary aperiodic uniformly recurrent word which is not
Sturmian, admits infinitely many minimal subshifts in its abelian
closure. Moreover, in the case of rational letter frequency, the
abelian closure  always contains a morphic image of the full shift.

In general, the abelian closure of an infinite word might have a
pretty complicated structure. T. Hejda, W. Steiner, and L.Q. Zamboni
studied the abelian shift of the Tribonacci word $\infw{T}$ defined
as the fixed point of $\tau:0\mapsto 01$, $1\mapsto 02$, $2\mapsto 0$.
They have announced that $\abclsr{\infw{T}}$ contains only one minimal
subshift, namely $\soc{\infw{T}}$ itself, but that there exist other
words in it as well \cite{HejdaSteinerZamboni15,ZamboniPersonal}.

The study of abelian closures is motivated by the question of, given
an infinite word $\infw{x}$, how strong is the bond between its abelian factors
and its language. We quantify this bond by the size of the abelian closure.
By size we do not mean the usual cardinality of a subshift, rather, we mean
the number of disjoint \emph{minimal} subshifts contained in $\abclsr{\infw x}$.
A shift orbit closure $\soc{\infw{x}}$ is minimal if does not properly contain
another shift orbit closure. If $\abclsr{\infw x}$ is huge (it contains
infinitely many minimal subshifts), then this bond is quite weak. On the other
hand, if $\abclsr{\infw x}$ is small (finitely many minimal subshifts), then
the bond is quite strong. The strongest bond is attained when $\abclsr{\infw{x}}$
is a minimal subshift itself. In this case we necessarily have
$\abclsr{\infw{x}} = \soc{\infw{x}}$. It is not hard to see that for purely
periodic words, their abelian closure is finite (see \autoref{prop:periodicFinite}).
On the other hand, the abelian closure of an ultimately periodic word can be huge
(see \autoref{ex:periodicHuge}). For reasons stemming from this observation, when
dealing with abelian closures of individual words, we shall assume the words to
define minimal shift orbit closures.

There is indeed no particular reason to restrict the definition of abelian closures
to just individual words; the abelian closure of a set of words $X$ comprises those
infinite words $\infw{y}$ whose each factor is abelian equivalent to some factor of
one of the words in $X$. For example, the abelian closure of $\soc{\infw{x}}$ coincides
with $\abclsr{\infw{x}}$. As mentioned previously, the shift orbit closure of an
infinite word is a certain type of a \emph{subshift}. In general, a non-empty set
$X \subseteq \Sigma^{\N}$ of infinite words is called a subshift if it is closed
(as a subset of the compact metric space $\Sigma^{\N}$ equipped with the usual product
topology defined by the discrete topology on $\Sigma$) and that $\sigma(X) \subseteq X$,
where the shift map $\sigma$, is defined as $\sigma(\infw{x})_i = \infw{x}_{i+1}$.%
{\footnote{Subshifts are often defined over bi-infinite words, in which case we require $\sigma(X) = X$ in the definition.}}
A subshift is called \emph{minimal} if it does not contain a proper subshift. Hence
a minimal subshift is always the shift orbit closure of some word $\infw{x}$.
It is routinely checked that $\abclsr{X}$ is a subshift for any set $X$ of words.

In this paper, we study how the characterization of Sturmian words as aperiodic
uniformly recurrent words with $\abclsr{\infw x} = \soc{\infw x}$ extends
to non-binary alphabets. We then study the abelian closures of certain
generalizations of Sturmian words, and some preliminary results have been
reported at DLT 2018 conference \cite{DBLP:conf/dlt/KarhumakiPW18}. Besides that,
we discuss abelian closures of subshifts in general. In \autoref{sec:bal} 
we characterize the abelian subshifts of aperiodic recurrent balanced words; they
are a finite union of minimal subshifts. In \autoref{sec:mincompl}  we
consider abelian closures of words over a $k$-letter alphabet with factor complexity
$n+k-1$ for each $n$, which are aperiodic words of minimal complexity involving $k$
letters. The behavior is different depending on $k$. For $k=2$, we are in the case
of Sturmian words, so we have  $\abclsr{\infw x}=\soc{\infw x}$. Surprisingly, the
most complicated behaviour is exhibited in the ternary alphabet. We show that for $k=3$, depending
on the word $\infw x$, its abelian closure $\abclsr{\infw x}$ contains either exactly
one, or uncountably many minimal subshifts. For alphabets of size greater than $3$,
$\abclsr{\infw x}$ equals the union of exactly two minimal subshifts, $\soc{\infw x}$
and its "reversal". Further, in \autoref{sec:AR}, we show that for Arnoux--Rauzy words,
their abelian closures contain non-recurrent words, and hence
$\abclsr{\infw x} \neq \soc{\infw x}$. We then extend our interest to general subshifts
in \autoref{sec:general_subshifts}. Our focus is on subshifts defined using notions from formal language
theory. We show that the abelian closure of a \emph{subshift of finite type} (resp., \emph{sofic subshift})
is not necessarily a subshift of finite type (resp., a sofic subshift) (see \autoref{sec:general_subshifts} for definitions).
We then conclude with open problems.

\section{Notation and first observations}
We recall some notation and basic terminology from the literature of combinatorics on
words. We refer the reader to \cite{lothaire1983combinatorics,MR1905123} for more on the subject.

The set of finite words over an alphabet $\Sigma$ is denoted by
$\words$ and the set of non-empty words is denoted by $\Sigma^+$.
The empty word is denoted by $\eps$. 
We let $|w|$ denote the length of a word
$w\in\words$. By convention, $|\eps| = 0$. 
A \emph{factor} of a word $\infw x$ is any block of its
consecutive letters, and we let $\mathcal L({\infw x})$ denote the
set of factors of $\infw x$. The length $n$ factors of $\infw x$ is
denoted by $\mathcal L_n({\infw x})$. The length $n$ prefix of the word $\infw{x}$ is denoted by
$\pref_n(\infw{x})$. The \emph{factor complexity} function $\complfunction{\infw{x}}$
is defined by $\complfunction{\infw{x}}(n) = | \mathcal L_n(\infw x)|$.
An infinite word $\infw x$ is
called \emph{recurrent} if each factor of $\infw x$ occurs infinitely
many times in $\infw x$. Further, $\infw x$ is \emph{uniformly recurrent}
if for each factor $u\in \mathcal L(\infw x)$ there exists $N\in \N$ such
that $u$ occurs as a factor in each factor of length $N$ of $\infw x$.
For a finite word $u\in\words$, we let $|u|_a$ denote the number of
occurrences of the letter $a\in\Sigma$ in $u$. For a finite word $v$, we let 
$v^\omega$ denote the infinite word obtained by repeating $v$ infinitely many 
times.

For $\infw x\in \Sigma^{\N}$ and $a\in\Sigma$, the limits
\begin{equation*}
\supfreq[\infw{x}]{a} := \lim_{n\to\infty}\frac{\sup_{v\in \mathcal L_n(\infw x)}|v|_a}{n}
\quad
\text{and}
\quad
\inffreq[\infw{x}]{a} := \lim_{n\to\infty}\frac{\inf_{v\in \mathcal L_n(\infw x)}|v|_a}{n}
\end{equation*}
exist.
and, moreover,
\begin{equation*}
\supfreq[\infw{x}]{a} = \inf_{n\in\N}\frac{\sup_{v\in \mathcal L_n(\infw x)}|v|_a}{n} \quad \text{and} \quad
\inffreq[\infw{x}]{a} = \sup_{n\in\N}\tfrac{\inf_{v\in \mathcal L_n(\infw x)}|v|_a}{n}.
\end{equation*}
This follows from Fekete's lemma as $\sup_{v\in\mathcal L_n(\infw x)}|v|_a$
(resp., $\inf_{v\in\mathcal L_n(\infw x)}|v|_a$) is \emph{subadditive} (resp., \emph{superadditive}) with respect to $n$.
It immediately follows that
\begin{equation}\label{eq:finiteLengthFrequency}
\sup_{v\in \mathcal L_n(\infw x)}|v|_a \geq n\, \supfreq[\infw{x}]{a}
 \quad \text{and} \quad
 \inf_{v\in \mathcal L_n(\infw x)}|v|_a \leq n\, \inffreq[\infw{x}]{a}
\end{equation}
 for all $n\in \N$.
If $\supfreq[\infw{x}]{a} = \inffreq[\infw{x}]{a}$, we denote the common limit by $\freq[\infw x]{a}$ and we say that \emph{$\infw x$ has uniform frequency of $a$}.

A \emph{subshift} $X\subseteq \Sigma^{\N}$, $X\neq \emptyset$, is
a closed set (with respect to the product topology on
$\Sigma^{\N}$) satisfying $\sigma(X)\subseteq X$, where $\sigma$
is the \emph{shift operator} (defined by $\sigma(a_0a_1a_2\cdots)
= a_1a_2\cdots$). For a subshift $X\subseteq \Sigma^{\N}$ we let
$\mathcal L(X) = \cup_{\infw y\in X}\mathcal L(\infw y)$. A subshift
$X\subseteq \Sigma^{\N}$ is called \emph{minimal} if $X$ does not
contain any proper subshifts. Observe that two minimal subshifts
$X$ and $Y$ are either equal or disjoint. Let $\infw x \in
\Sigma^{\N}$. We let $\soc{\infw x}$ denote the \emph{shift
orbit closure} of $\infw x$, which may be defined as the subshift
$\{\infw y\in \Sigma^{\N} \colon \mathcal L(\infw y)\subseteq
\mathcal L(\infw x)\}$. Thus $\mathcal L(\soc{\infw x}) =
\mathcal L(\infw x)$ for any word $\infw x\in\Sigma^{\N}$. It is
known that $\soc{\infw x}$ is minimal if and only if $\infw x$ is
uniformly recurrent. For a morphism $\varphi:\Sigma\to \Delta^*$
(that is, $\varphi(uv) = \varphi(u)\varphi(v)$ for all
$u,v\in\words$) and a subshift $X\subseteq \Sigma^{\N}$, we define
$\varphi(X) = \cup_{\infw x\in X}\soc{\varphi(\infw x)}$. When
using \emph{erasing morphisms}, that is, some letter maps to
$\eps$, we make sure that no point in $X$ gets mapped to a finite
word. For more on this topic we refer the reader to
\cite{LindMarcus95}.

We recall definitions and properties related to
\emph{Sturmian words} from \cite[Chapter 2]{MR1905123}.
We identify the interval $[0,1)$ with the unit circle $\mathbb T$
(the point $1$ is identified with point $0$). For points
$x,y\in\mathbb T$, we let $\underline{I} (x,y)$ (resp.,
$\bar{I}(x,y)$) denote the half-open interval on $\mathbb T$
containing $x$ (resp., $y$) and starting from $x$ and ending at
$y$ in counter--clockwise direction. We omit the bar whenever it
does not matter which endpoint is in the interval. Let $\alpha\in
\mathbb{T}$ be irrational or rational and let $\rho\in \mathbb T$.
The map $R_{\alpha}:\mathbb T \to \mathbb T$, $x \mapsto
\{x+\alpha\}$, where $\{x\} = x-\lfloor x \rfloor$ is the
fractional part of $x\in\R$, defines a (counter-clockwise) rotation on $\mathbb T$.
Divide $\mathbb T$ into two half-open intervals $I_0 =
\underline{I}(0,1-\alpha)$ and $I_1 =
\underline{I}(1-\alpha,1)$ (resp., $I_0 = \bar
I(0,1-\alpha)$ and $I_1 = \bar I(1-\alpha,1)$) and define the
coding $\nu:\mathbb T\to \{0,1\}$, $x\mapsto i$ if $x\in I_i$,
$i=0,1$. The \emph{rotation word} $\underline{\infw s}_{\alpha,\rho}$
(resp., $\bar{\infw s}_{\alpha,\rho}$) of \emph{slope} $\alpha$ and
\emph{intercept} $\rho$ is the word $a_0a_1\cdots \in
\{0,1\}^{\N}$ defined by $a_n = \nu(R_{\alpha}^n(\rho))$ for all
$n\in\N$.
Note that $00$ occurs in $\infw s_{\alpha,\rho}$ if
and only if $\alpha<1/2$. Clearly, $\infw s_{\alpha,\rho}$ is aperiodic if
and only if $\alpha$ is irrational. Each aperiodic rotation word is a Sturmian
word and vice versa. We call periodic rotation words periodic Sturmian. Observe the special role played by
$\infw s_{\alpha,\alpha}=\infw s$:
both $01\infw s$ and $10\infw s\in \soc{\infw s}$ for $\alpha\in(0,1)$.

Every Sturmian word $\infw s$ is uniformly recurrent so that
$\soc{\infw s}$ is minimal. Further, $\infw s'\in \soc{\infw s}$ if
and only if $\infw s$ and $\infw s'$ are of the same slope. In particular,
the intercepts or the endpoints of $I_0$ and $I_1$ do not play any
role when speaking of the shift orbit closure of a Sturmian word.
(In Section \ref{sec:mincompl} we pay attention to these choices.)

An infinite word $\infw x\in\Sigma^{\N}$ is called \emph{balanced} if, for all
$v,v'\in \Sigma^{\N}$ with $|v| = |v'|$ and for all $a\in \Sigma$, we have
$||v|_a-|v'|_a|\leq 1$.
Periodic and aperiodic Sturmian words are exactly the recurrent balanced binary words
\cite{Morse10.2307/2371431}. It follows that, for each (aperiodic or periodic) Sturmian
word $\infw s$, the set $\{|v|_1 \colon v\in \mathcal L_n(\infw s)\}$ consists of at most
two values $k$ and $k+1$ for some $k$ depending on $\infw s$ and $n$. Observe
now that $\freq[\infw s]{1} = \alpha$ for any Sturmian word $\infw s$ of slope $\alpha$. By \eqref{eq:finiteLengthFrequency}, $k$ above equals $\lfloor n \alpha \rfloor$.

We turn to the main notion of this paper. 

\begin{definition} For $\infw x\in
\Sigma^{\N}$ we define the \emph{abelian closure} of $\infw x$ as
$$\abclsr{\infw x}= \{\infw y\in \Sigma^{\N}\mid \forall u\in
\mathcal L(\infw y) \exists v\in \mathcal L(\infw x)\colon
u\sim_{\text{ab}}v\}.$$\end{definition} 

Now, for any $\infw x\in
\Sigma^{\N}$, the abelian closure $\abclsr{\infw x}$ is indeed
a subshift.
 We make preliminary observations on
abelian closures of infinite words.

\begin{lemma}\label{lem:uniformFrequencies}
Assume $\infw x\in\Sigma^{\N}$ has uniform frequency of a letter
$a\in\Sigma$. Then any word $\infw y\in\abclsr{\infw x}$ has
uniform frequency of $a$ and $\freq[\infw y]{a} = \freq[\infw x]{a}$.
\end{lemma}
\begin{proof}
For all $\infw y\in\abclsr{\infw x}$ and for all $n\in\N$, we
have immediately from the definition of $\abclsr{\infw x}$ that
\begin{equation*}
\sup_{v\in \mathcal L_n(\infw x)}\frac{|v|_a}{n}
\geq \sup_{v\in \mathcal L_n(\infw y)}\frac{|v|_a}{n}
\geq \inf_{v\in \mathcal L_n(\infw y)}\frac{|v|_a}{n}
\geq \inf_{v\in \mathcal L_n(\infw x)}\frac{|v|_a}{n}.              
\end{equation*}
Letting $n\to\infty$ gives our claims.
\end{proof}
We immediately have that if $\infw x$ has an irrational uniform frequency
of some letter $a$, then $\abclsr{\infw x}$ contains only aperiodic words.
We continue by observing how the abelian closures of
periodic and ultimately periodic words can differ.

\begin{proposition}\label{prop:periodicFinite}
For any periodic word $\infw x$, the abelian closure $\mathcal
A({\infw x})$ is finite.
\end{proposition}
\begin{proof}
A word $\infw y$ is periodic if and only if all factors of length $n$ are abelian equivalent
for some $n\geq 1$ \cite{DBLP:journals/mst/CovenH73}. Let $n$ be the least such integer
for $\infw x$. It follows that all factors of length $n$ of any word $\infw y \in \mathcal
A({\infw x})$ are abelian equivalent. Thus $\infw{y} = v^{\omega}$ with $|v|$ dividing $n$. There
are finitely many such words.
\end{proof}
In general, the abelian closure of an ultimately periodic word can
be huge.
\begin{example}\label{ex:periodicHuge}
Let $\infw x = 0011(001101)^{\omega}$.
It is readily verified that for odd lengths $\infw x$ has two abelian factors, and
for even lengths three. Further, for each factor of $\infw x$, the number of occurrences of $1$
differs by at most one from half of its length.
Thus, by the discussion in the introduction, we have
$\infw{TM}\in \abclsr{\infw x}$ so that $\abclsr{\infw x} = \abclsr{\infw{TM}} = \{\eps,0,1\}\cdot\{01,10\}^{\N}$.
\end{example}
A family of such examples are given in \cite[Ex.~2]{DBLP:conf/dlt/KarhumakiPW18}.

In the end of this section we show that, for a uniformly recurrent binary word $\infw x$,
$\abclsr{\infw x}$ contains exactly one minimal subshift
if and only if $\infw x$ is a Sturmian word.
We start with a crucial observation, which is characteristic only for binary words.
\begin{lemma}[Corridor Lemma]\label{lem:corridorLemma}
Let $\infw x$ be a binary word. Then $\infw y\in \abclsr{\infw x}$ if and only if, for
all $n\in\N$, we have
\begin{equation*}
\inf_{u \in \lang[n]{\infw {y}}}|u|_1 \geq \inf_{u \in \lang[n]{\infw {x}}}|u|_1
\quad
\text{and}
\quad
\sup_{u \in \lang[n]{\infw {y}}}|u|_1 \leq \sup_{u \in \lang[n]{\infw {x}}}|u|_1.
\end{equation*}
\end{lemma}
\begin{proof}
It is easy to see (e.g., by a sliding window argument) that, for any $n\geq 1$,
there exists a word $u\in \mathcal L_n(\infw z)$ with $|u|_1 = m$ if and only if
$\inf_{v \in \lang[n]{\infw{x}}} |v|_1 \leq m \leq \sup_{v \in \lang[n]{\infw{x}}}|v|_1$.
Applying this observation to $\infw x$ and $\infw y$ we have that for each
$v\in \mathcal L_n(\infw y)$ there exists $u\in \mathcal L_n(\infw x)$ such that
$v\sim_{\text{ab}} u$ if and only if
$\inf_{u \in \lang[n]{\infw {x}}}|u|_1  \leq \inf_{u \in \lang[n]{\infw {y}}}|u|_1$
and
$\sup_{u \in \lang[n]{\infw {y}}}|u|_1\leq \sup_{u \in \lang[n]{\infw {x}}}|u|_1$.
\end{proof}

We now characterize Sturmian words in terms of abelian closures. 
\begin{theorem}\label{th:St}
Let $\infw x\in\{0,1\}^{\N}$ be uniformly recurrent. Then $\mathcal
A({\infw x})$ contains exactly one minimal subshift if and only if $\infw x$ is Sturmian.
\end{theorem}

Another version of this is the following:

\begin{theorem}\label{th:St1}
Let $\infw x\in\{0,1\}^{\N}$ be aperiodic and uniformly recurrent. Then $\mathcal
A({\infw x})=\soc{\infw x}$ if and only if $\infw x$ is Sturmian.
\end{theorem}

\begin{proof}
First we show that $\abclsr{\infw x} = \soc{\infw x}$ for
Sturmian words. Let $\infw x$ be a Sturmian word of slope $\alpha$
and $\infw y\in \abclsr{\infw x}$. By the
\hyperref[lem:corridorLemma]{Corridor Lemma}, $\infw y$ is balanced
and, by \autoref{lem:uniformFrequencies}, has uniform frequencies
of letters equal to those of $\infw x$. Thus $\infw y\in
\soc{\infw x}$.

Assume then that $\abclsr{\infw x}$ contains exactly one minimal subshift, namely $\soc{\infw{x}}$, and let $\alpha = \supfreq[\infw{x}]{1}$.
Take a (periodic or aperiodic) Sturmian word $\infw s$ of slope
$\alpha$. By the \hyperref[lem:corridorLemma]{Corridor Lemma},
we have $\soc{\infw s}\subseteq \abclsr{\infw x}$ using \eqref{eq:finiteLengthFrequency}.
Since $\soc{\infw s}$ is also minimal, we have $\soc{\infw x}=\soc{\infw s}$.
Thus $\infw x$ is Sturmian, and $\abclsr{\infw{x}} = \soc{\infw{x}}$.
\end{proof}

The following example shows that we cannot omit the assumption of uniform recurrence
from the statement of the above theorem.
\begin{example}
Take the \emph{Champernowne word} $\infw{C}_2$ (over the binary alphabet), which is obtained
by concatenating all finite words ordered by length and lexicographic order for the same length:
\begin{equation*}
\infw{C}_2= 0 \gap 1 \gap 00 \gap  01 \gap 10 \gap 11 \gap  000 \gap 001 \gap 010
\gap  011 \gap 100 \gap 101 \gap  110 \gap 111 \cdots
\end{equation*}
Clearly, both $\soc{\infw{C}_2}$ and $\abclsr{\infw{C}_2}$ are equal to the full shift, i.e., contain all binary words.
\end{example}

Note that the property $\abclsr{\infw x}=\soc{\infw x}$ or $\abclsr{\infw x}$ containing exactly
one minimal subshift does not characterize Sturmian words among
uniformly recurrent words over arbitrary alphabets. Let $\infw{f}$ be the \emph{Fibonacci word}, which is a Sturmian word defined as the fixed point of the morphism $0 \mapsto 01$,
$1\mapsto 0$. 
Let then $\varphi: 0\mapsto 02, 1\mapsto 12$. Then $\abclsr{\varphi(\infw f)} = \soc{\varphi(\infw f)}$ (see  \autoref{thm:abelianSSMinimalComplexity}).

We investigate possible generalizations of the property $\mathcal
A({\infw x})=\soc{\infw x}$ to nonbinary alphabets in the next sections. 

\section{Abelian closures of balanced words}\label{sec:bal}

In this section we study the abelian closures of non-binary aperiodic balanced words. We
prove that the abelian closure of such a word is a finite union of minimal subshifts:

\begin{theorem}\label{thm:BalancedabelianSS}
Let $\infw u$ be aperiodic recurrent and balanced. Then $\abclsr{\infw u}$ is the union of
finitely many minimal subshifts.
\end{theorem}
As we will show below, abelian closure of a recurrent balanced word can contain one or more (yet a finite number) of minimal subshifts, depending on its structure, and we can in fact compute this number. Our results rely heavily on the characterization of aperiodic recurrent balanced words by R.~Graham \cite{DBLP:journals/ipl/Graham72} and
P.~Hubert \cite{Hubert00}. In fact, the characterization allows us to characterize the
abelian closures of slightly more general words. In particular, the the techniques used in the proof of
the above theorem give us, for each $k$, 
an aperiodic word $\infw x_k$ over a four-letter alphabet
such that $\abclsr{\infw x_k}$ equals the union of $k$ distinct minimal subshifts (see Proposition \ref{prop:kMinimalSubshifts}).

We need some notation to give a characterization of aperiodic recurrent words.

\begin{definition}
A word is called \emph{constant gap} if each letter occurs
with a constant gap.
\end{definition}
For example, $(abac)^{\infty}$ is a constant gap word.

\begin{definition}
Let $\infw x$ be a finite or infinite binary word, $\infw z_0\in A^{\N}$ and
$\infw z_1\in B^{\N}$, where $A$ and $B$ are some alphabets.
Let $\mathcal S(\infw x,\infw z_0,\infw z_1)$ denote the word obtained
from $\infw x$ by substituting the $n$th occurrence of $0$ (resp., $1$) in $\infw{x}$ by the $n$th letter of $\infw z_0$ (resp., $\infw{z}_1$).
\end{definition}
We illustrate the above operation with an example.
\begin{example}
Let $\infw f=01001010010010100101\cdots$ be the Fibonacci word,
$\infw z_0=(0102)^{\omega}$, and $\infw z_1=(ab)^{\omega}$. Then
$\mathcal S(\infw f,\infw z_0,\infw z_1)= 0a10b2a01b02a0b10a2b\cdots$
is balanced.
\end{example}

The following theorem characterises recurrent balanced words using constant gap words
and the operation $\mathcal{S}$.

\begin{theorem}[{\cite[Thm.~1]{Hubert00}}, \cite{DBLP:journals/ipl/Graham72}]\label{thm:characterizationBalance}
An aperiodic word $\infw u\in \Sigma^{\N}$ is recurrent and balanced if and only if there
exist a partition $\{A,B\}$ of $\Sigma$, two constant gap words $\infw z_0\in A^{\N}$ and
$\infw z_1 \in B^{\N}$, and a Sturmian word $\infw s$, such that
$\infw u = \mathcal S(\infw s,\infw z_0,\infw z_1)$.
\end{theorem}

We remark that although the structure of aperiodic balanced words is
clear, the structure of periodic balanced words is a mystery: the
following conjecture by Fraenkel, 1973, remains open despite efforts
of different scientists: The unique (up to a permutation of letters)
balanced word on $k\geq 3$ letters with all distinct frequencies of
letters is $(F_k)^{\omega} = (F_{k-1}kF_{k-1})^{\omega}$ where
$F_2 = 121$ \cite{Fraenkel73}. The conjecture has been verified for
$k\leq 7$ (see \cite{BARK03} and references therein).

In fact, throughout this section we consider slightly more general words, namely, we relax
the condition of $\infw z_0$ and $\infw z_1$ being constant gap. We obtain a
characterization of the abelian closures of such words, from which the characterization
of abelian closures of recurrent balanced words follows.

We need the following lemma:

\begin{lemma}
Let $\infw u = \mathcal S(\infw s,\infw z_0,\infw z_1)$ with $\infw s$
Sturmian and $\infw z_i$ periodic words. Then
\begin{equation*}
\mathcal L(\infw u) = \{\mathcal S(x,\infw z_0',\infw z_1')\colon x\in \mathcal L(\infw s), \infw z_i'\in \soc{\infw z_i}\}.
\end{equation*}
Further, $\infw u$ is uniformly recurrent.
\end{lemma}
\begin{proof}
In fact, this result is implicitly contained in \cite{Hubert00}, only
stated in slightly weaker form. Theorem 2 and Proposition 3.1 of \cite{Hubert00}
essentially state the following. Let
$\infw{v} = \mathcal{S}(\infw{s},\infw{y}_0,\infw{y}_1)$, where $\infw{s}$ is 
a Sturmian word and $\infw y_0$ and $\infw{y}_1$ are constant gap words with 
periods $l_0$ and $l_1$ respectively. Then the factor complexity function of 
$\infw v$ satisfies $\complfunction{\infw v}(n) = l_0l_1(n+1)$ for all large 
enough $n$. Further \cite[Prop.~5.1]{Hubert00} states that $\infw v$ is
uniformly recurrent.

Observe that
$\mathcal L(\infw v) \subseteq
	\{\mathcal S(x,\sigma^i(\infw y_0),\sigma^j(\infw y_1))\colon x\in \mathcal L(\infw s), i,j\in\N\}$
and the cardinality of the latter set equals $l_0l_1(n+1)$ for large enough
$n$. This coincides with the factor complexity $\complfunction{\infw v}(n)$ by
Hubert's result.
We deduce that for all $x\in \mathcal L(\infw s)$, $i,j\in\N$, we have
$\mathcal S(x,\sigma^i(\infw y_0),\sigma^j(\infw y_1)) \in \mathcal L(\infw v)$.

Now we can take the constant gap words $\infw y_0 = (a_1\dots a_{l_0})^{\omega}$ and
$\infw y_1 = (b_1\dots b_{l_1})^{\omega}$, where $l_0$ and $l_1$ are the lengths of
the periods of $\infw z_0$ and $\infw z_1$: $\infw z_i = u_i^{\omega}$, $|u_i|=l_i$.
Define a coding $\tau$ so that $a_1\cdots a_{l_0} \mapsto u_0$ and $b_1\cdots b_{l_1}\mapsto u_1$.
Then $\tau(\infw v)=\infw u$ so that $\infw u$ is uniformly recurrent. Further,
$\tau (\mathcal S(x,\sigma^i(\infw y_0),\sigma^j(\infw y_1)))=\mathcal S(x,\sigma^i(\infw z_0),\sigma^j(\infw z_1))$. The claim follows now immediately.
\end{proof}

In the above lemma we allow the periodic words $\infw z_0$ and $\infw z_1$ contain common letters.
On the other hand, in the following proposition we assume that they do not share common letters.
This puts us in the position of characterizing the abelian closures of such words.

\begin{proposition}\label{prop:abSSgeneralizedBalanced}
Let $\infw u=\mathcal S(\infw s,\infw z_0,\infw z_1)$ for some Sturmian word $\infw s$ and periodic
words $\infw z_0\in A^{\N}$ and $\infw z_1\in B^{\N}$, where $A$ and $B$ are disjoint alphabets. 
Then
\begin{equation*}
    \abclsr{\infw u} = \bigcup_{\infw t_i\in \abclsr{\infw z_i}}
														\soc{\mathcal S(\infw s,\infw t_0,\infw t_1)},
\end{equation*}
and $\abclsr{\infw u}$ is
a finite union of minimal subshifts.
\end{proposition}
\begin{proof}
Let first $\infw x\in \soc{\mathcal S(\infw s,\infw t_0,\infw t_1)}$ for some
$\infw t_i\in \abclsr{\infw z_i}$.
For any factor $x\in \mathcal L(\infw x)$, we have
$x = \mathcal S(y,\sigma^{i}(\infw t_0),\sigma^j(\infw t_1))$ for some
$y\in \mathcal L(\infw s)$ and $i,j\geq 0$ by the above lemma. Let $u = \pref_n(\sigma^i(\infw t_0))$
and $v = \pref_m(\sigma^j(\infw t_1))$, where $n = |y|_0$, $m=|y|_1$. By assumption,
there exist $u'\in \mathcal L(\infw z_0)$ and $v'\in \mathcal L(\infw z_1)$ such that
$u\sim_{\text{ab}} u'$ and $v\sim_{\text{ab}} v'$. Choose $r,s$ such that
$\sigma^r(\infw z_0)$ begins with $u'$ and $\sigma^s(\infw z_1)$ begins with $v'$. It
follows that $x\sim_{\text{ab}} \mathcal S(y,\sigma^r(\infw z_0),\sigma^s(\infw z_1))\in \mathcal L(\infw u)$.
We thus have $\infw x\in \abclsr{\infw u}$.

Let then $\infw x\in \abclsr{\infw u}$. Take $\varphi:\Sigma\to \{0,1\}$ such that
$\varphi(a) = 0$ if and only if $a\in A$. It follows that
$\varphi(\infw x)\in \abclsr{\varphi(\infw u)} = \soc{\infw s}$ since the alphabets
$A$ and $B$ are disjoint.
Take then the morphism $\varphi_A:\Sigma \to A^*$ such that $\varphi_A(a) = a$ for $a\in A$,
otherwise $\varphi_A(a) = \eps$. Define the morphism $\varphi_B \colon \Sigma \to B^*$ analogously. We 
again have that $\varphi_A(\infw x) \in \abclsr{\varphi_A(\infw u)}$
where $\varphi_A(\infw u) = \infw z_0$. Similarly $\varphi_B(\infw x) \in \abclsr{\infw z_1}$.
It is now evident that $\infw x = \mathcal S(\infw s',\infw t_0,\infw t_1)$ for some
$\infw t_i \in \abclsr{\infw z_i}$, $i=0,1$, and $\infw s' \in \soc{\infw s}$.
The above lemma implies that $\infw x\in \soc{\mathcal S(\infw s,\infw t_0,\infw t_1)}$. 

As $\soc{\infw{z}_i}$ is finite by \autoref{prop:periodicFinite}, 
 $\abclsr{\infw u}$ is
a finite union of minimal subshifts.
This concludes the proof.
\end{proof}

The above proposition has \autoref{thm:BalancedabelianSS} as an immediate corollary.

Let us consider what the above proposition says. The number of distinct minimal subshifts is 
bounded above by the product of the number of minimal subshifts in $\abclsr{\infw{z}_0}$ and the number of minimal subshifts in $\abclsr{\infw{z}_1}$.

\begin{example}
Let $\infw z_0 = (0102)^{\omega}$ and $\infw z_1 = (34)^{\omega}$. Now $\abclsr{\infw z_i} = \soc{\infw z_i}$
as is readily verified. Thus $\abclsr{\infw u} = \soc{\infw u}$ for $\infw u = \mathcal S(\infw s,\infw z_1,\infw z_2)$,
$\infw s$ Sturmian, by the above proposition. 
\end{example}

For the case of recurrent balanced words, the periodic words in the construction are constant
gap words. It is natural to ask whether the extra property of constant gaps restricts the
cardinality of the number of minimal subshifts in its abelian closure. We give a negative answer to this:

\begin{example}
 Let $\infw z_1 = (a_0 a_1 a_2 a_3 a_4 A_0 \cdot a_0 a_1 a_2 a_3 a_4 A_1 \cdot a_0 a_1 a_2 a_3 a_4 A_2)^{\omega}$. Here the letters $a_i$ have constant gaps of length $6$, and $A_i$ have constant gaps of length $18$. Take any
constant gap sequence $\infw z_2$ which is not closed under reversal (e.g., $(abc)^{\omega}$), and $\infw u = \mathcal S(\infw s,\infw z_1,\infw z_2)$. Then we can independently take reversal inside $\infw z_1$, inside $\infw z_0$, and inside the arithmetic progression given by $A_0A_1A_2$ in $\infw z_0$. We thus get eight minimal subshifts. Note that this construction can be generalized to produce $2^k$ minimal subshifts for any $k$.
\end{example}

The techniques used in the proof of
the above theorem give us the following proposition. We remark that the words in question are
not necessarily balanced:

\begin{proposition}\label{prop:kMinimalSubshifts}
For each $k \geq 1$ there exists an aperiodic word $\infw x_k$ over a three-letter alphabet
such that $\abclsr{\infw x_k}$ equals the union of $k$ distinct minimal subshifts.
\end{proposition}

\begin{proof}
For $k=1$ we may take any Sturmian word. Let thus $k\geq 2$.
Consider first the abelian closure of the periodic word $\infw z = (0^{2k-1}11)^{\omega}$.
It is readily verified that $\abclsr{\infw x}$ contains the minimal subshifts generated
by the words $(0^{2k-1-i}10^{i}1)^{\omega}$, $i=0,\ldots,k-1$. To see that 
there is nothing else in $\abclsr{\infw z}$, we observe the following.
By \autoref{prop:periodicFinite}, any
word in $\abclsr{\infw{z}}$ is periodic with period dividing $2k+1$, an odd
number. The period cannot be less than $2k+1$, as the number of $1$s in
factors of length $2k+1$ is only $2$. On the other hand, all words of length
$2k+1$ containing two occurrences of $1$ (and hence the periodic words they
generate) occur already in the subshifts above.

For the claim we set $\infw x_k = \mathcal S(\infw s,\infw z,a^{\omega})$
where $\infw s$ is an aperiodic Sturmian word.
By \autoref{prop:abSSgeneralizedBalanced},
\begin{equation*}
\abclsr{\infw x_k} = \bigcup_{i=0}^{k-1}\soc{\mathcal S(\infw s,(0^{2k-1-i}10^{i}1)^{\omega},a^{\omega})},
\end{equation*}
a union of $k$ distinct minimal subshifts.
\end{proof}

\section{Abelian closures of words of minimal complexity}\label{sec:mincompl}

First we study the abelian closures of aperiodic nonbinary words of minimal factor
complexity. Over an alphabet $\Sigma$ with at least two letters, the minimal complexity is $n+|\Sigma|-1$. The
structure of words of complexity $n+C$ is related to the structure of Sturmian words and is
well understood (\cite{FM97,Didier99,KaboreTapsoba07}). 
The main goal of this subsection is to prove that for aperiodic ternary words of minimal
complexity their abelian closure consists of either one or uncountably many minimal 
subshifts (\autoref{thm:abelianSSMinimalComplexity}); for alphabets of size greater than $3$ 
the abelian closure contains exactly two minimal subshifts (when the words are assumed to be
recurrent, \autoref{th:min_compl}). Recall that for binary alphabet, we have exactly one
minimal subshift (\autoref{th:St}).

A proof of the following is explained in \cite{FM97} in the discussion following Lemma 1.
\begin{lemma}
A minimal complexity word $\infw{u}$ over the alphabet $A$ is of the form
$a_0\cdots a_t \infw{u}'$, where $\infw{u}'$ is a recurrent minimal complexity
word over an alphabet $A' \subseteq A$, and $a_0$, \ldots, $a_t$ are distinct
letters of $A \setminus A'$.
\end{lemma}

We shall consider the abelian subshifts of recurrent aperiodic minimal complexity words. The following lemma then extends those results to handle
non-recurrent ones as well.
\begin{lemma}\label{lem:abclsrNonRecurrent}
Let $\infw{u} = a_0\cdots a_t \infw{u}'$ be an aperiodic minimal complexity word as in the above lemma. Then $\abclsr{\infw{u}} = \soc{\infw{u}} \cup \abclsr{\infw{u'}}$.
\end{lemma}
\begin{proof}
Let $\infw{y} \in \abclsr{\infw{u}}$, so it contains at most one occurrence of
each of the letters $a_i$. If it contains none of them, there is nothing to
prove.

So let us write $\infw{y} = p a_i \infw{y}'$, where $\infw{y}'$ contains none
of the letters $a_j$, $j=1$, \ldots, $t$. First of all, $i = t$, as $a_i$ must
be followed by $a_{i+1}$ (or $a_{i-1}$, in which case we reach $a_0$ at some
point, after which we cannot add anything). We show that
$\infw{y}' = \infw{u}'$. Assume that $\infw{y}'\neq \infw{u}'$, so $\infw{y}$
contains the factor $a_t w a$, while $\infw{u}$ contains $a_t w b$, for some
$a\neq b$. This is impossible, as $a_twb$ is the only factor of length $|w|+2$
of $\infw{u}$ that contains only letters from $A \setminus A'$ (apart from $a_t$).
Hence $\infw{y} = p a_t \infw{u}'$.

Now either $p$ ends with $a_{t-1}$ or with a letter $a  \in A'$. In the latter case, we must have that $a$ is the first
letter of $\infw{u}'$, and so $\infw{y}$ contains $a a_t a$. It follows that
$\infw{u}'$ must begin with $aa$, so $\infw{y}$ contains $a a_t aa$ and hence
$\infw{u}'$ must begin with $aaa$. By continuing with this line of reasoning,
we see that $\infw{u}' = a^{\omega}$, which is contrary to the assumption that $\infw{u}$ is aperiodic. We deduce that $p$ must end with $a_{t-1}$. Now $p$
cannot contain any letter from $A'$ (it would be followed by a letter $a_i$
with $i<t$, a contradiction). The only option is that
$\infw{y} = a_j a_{j+1} \cdots a_{t}\infw{u}'$ for some $j \geq 0$, which suffices for the proof. 
\end{proof}

\subsection{Ternary minimal complexity words}
We start with infinite words for which $p(n) = n+2$ for all $n\geq
1$. Observe that this implies that we deal with ternary words.

In \cite{Cassaigne98sequenceswith}, J. Cassaigne characterizes words having factor complexity $n+C$ for all $n\geq n_0$,
$C$ a constant. Here we consider the case of $C=2$ and $n_0=1$.
We first recall their characterization, which can be deduced from
\cite{KaboreTapsoba07} (see also \cite{Didier99}).

\begin{theorem}\label{thm:characterizationNp2Complexity}
A word $\infw u\in \{0,1,2\}$ has factor complexity $\mathcal P(n) = n+2$ for all $n\geq 1$
if and only if $\infw u$ is of the form (up to permuting the letters)
\begin{enumerate}[topsep=0pt]
\item\label{item:item1LowComplexity} $\infw u = 2 \infw s$ for some Sturmian word
$\infw s\in\{0,1\}^{\N}$, or
\end{enumerate}
$\infw u \in \soc{\varphi(\infw s)}$, where $\infw s$ is a Sturmian word
and $\varphi$ is defined by
\begin{enumerate}[resume,topsep=0pt]
\item \label{item:item2LowComplexity} $0\mapsto 02$, $1\mapsto 12$;
\item \label{item:item3LowComplexity} $0\mapsto 0$, $1\mapsto 12$.
\end{enumerate}
\end{theorem}
In this subsection we study the abelian closures of these words. The main result is the following theorem.

\begin{theorem}\label{thm:abelianSSMinimalComplexity}
Let $\infw u$ be a word of factor complexity $n+2$ for all $n\geq 1$. If
$\infw u$ is as in \autoref{thm:characterizationNp2Complexity} \eqref{item:item1LowComplexity} or \eqref{item:item2LowComplexity},
then $\abclsr{\infw u} = \soc{\infw u}$. If $\infw u$ is as in
\eqref{item:item3LowComplexity}, then $\abclsr{\infw u}$ contains uncountably many
minimal subshifts.
\end{theorem}

In fact we are able to characterize the abelian closures of these words.
We do this in parts, the first two cases are straightforward and we prove
them first. For the last case we need some further notions.
\begin{remark}
In the following, 
we often use the following argument. Assume that each letter in $\infw x\in \Sigma^{\N}$ occurs
with bounded gaps. Let $\varphi$ be a morphism such that $|\varphi(a)|\leq 1$. Then
$\varphi(\abclsr{\infw x})\subseteq \abclsr{\varphi(\infw x)}$. Indeed, since in $\infw x$ each letter occurs with bounded gaps, the same holds for any
$\infw y\in \abclsr{\infw x}$. Consequently $\varphi(\infw y)$ is infinite.
Letting $v$ be a factor of $\varphi(\infw{y})$, there exists a factor $v'$ of
$\infw{y}$ such that $\varphi(v') = v$ due to the length assumption on $\varphi$. As $\infw{y} \in \abclsr{\infw{x}}$ there exists a factor $w'$
of $\infw{x}$ abelian equivalent to $v'$. It thus follows that $\varphi(\infw{x})$ contains the factor $\varphi(w')$ abelian equivalent to
$v$. As $\infw{y}$ and $v$ were arbitrary, we conclude that $\varphi(\infw{y}) \in \abclsr{\varphi(\infw{x})}$.

In particular, if $\varphi(\infw x)$ is Sturmian, then $\varphi(\abclsr{\infw x})\subseteq \soc{\varphi(\infw x)}$
by \autoref{th:St}.
\end{remark}

\begin{proposition}\label{prop:characterizationNp2Small}
Let $\infw u$ as in \autoref{thm:characterizationNp2Complexity} \eqref{item:item1LowComplexity} or \eqref{item:item2LowComplexity}. Then $\abclsr{\infw u} = \soc{\infw u}$.
\end{proposition}
\begin{proof}Assume first that $\infw u = 2\infw s$. The claim follows immediately from \autoref{lem:abclsrNonRecurrent} together with \autoref{th:St1}.

Assume then that $\infw u $ is as in \eqref{item:item2LowComplexity}.
Let $\infw x\in \abclsr{\infw u}$. By applying the morphism $2\mapsto 2$, $1\mapsto0$, and
$0\mapsto 0$, we see that every second letter of $\infw x$ is $2$.
Further, by mapping $2\mapsto \eps$, $i\mapsto i$ for $i=0,1$,
we see that $\infw x$ maps into a word in $\soc{\infw s}$. It is straightforward to see that
now $\infw x\in \soc{\infw u}$.
\end{proof}

The rest of this subsection is devoted to the case where $\infw u$ is as in \autoref{thm:characterizationNp2Complexity} \eqref{item:item3LowComplexity}.
This case is more intricate as shown by the following example: $\abclsr{\infw u}$
contains non-recurrent words, similar to the Tribonacci word.
\begin{example}\label{ex:nonRecurrentInAbShift}
Let $\alpha \in \mathbb T$, let $\varphi$ be as in \autoref{thm:characterizationNp2Complexity} \eqref{item:item3LowComplexity}, and
$\infw u = \varphi(\infw s_{\alpha,\alpha})$. The words
$\infw u_1 = \varphi(01\infw s_{\alpha,\alpha}) = 012\infw u$ and
$\infw u_2 = \varphi(10\infw s_{\alpha,\alpha})= 120\infw u$ are
both in $\soc{\infw u}$. We claim that the non-recurrent word $\infw x = 02\infw u\in \abclsr{\infw u}$.
Indeed, $\sigma(\infw x) = \sigma^2(\infw u_1)\in\soc{\infw u}$ (recall $\sigma$ is the shift map). Further, any prefix of $\infw x$ of length
at least $2$ is abelian equivalent to the prefix of $\sigma(\infw u_2)$. Thus $\infw x\in \abclsr{\infw u}$.
\end{example}

We now analyze the structure of $\abclsr{\infw u}$. Without loss of
generality we may take $\infw u = \varphi(\infw s)$, since $\varphi(\infw s)$ is uniformly
recurrent. Consider the images of $\infw u$ under the morphisms $\varphi_1: 0\mapsto 0$,
$1\mapsto 1$, $2\mapsto 0$ and $\varphi_2:0\mapsto 0$, $1\mapsto 0$, $2\mapsto 1$. We have
$\varphi_1(\infw u) = \infw s_1$, a Sturmian word of some slope $\alpha$ and intercept $\rho$.
(Indeed, $\varphi(\infw u) = G \circ E (\infw s)$ using the notation of \cite[Chap.~2, p. 72]{MR1905123}.)
Symmetrically, $\varphi_2(\infw u) = \infw s_2$ is a Sturmian word
of slope $\alpha$ and intercept $\rho'$. (Again, $\varphi_2(\infw u) = D \circ E(\infw s)$, see again the above reference.) In
fact we can say more: $\rho' = \rho-\alpha$ or, equivalently, $\infw s_1 = \sigma(\infw s_2)$.
Observe now that $\infw s_1$ contains the factor $00$ meaning that $\alpha < 1/2$.

By mapping any word $\infw x\in \abclsr{\infw u}$ with the morphisms $\varphi_1$ and
$\varphi_2$, we obtain two Sturmian words with the same slope $\alpha$. Further, by applying
$0\mapsto \eps$ on $\infw u$, we see that $\infw x\in \soc{\infw (12)^{\omega}}$. This
implies that all words in $\abclsr{\infw u}$ are obtained by somehow "interleaving" two
Sturmian words of the same slope ($\infw s_2$ is encoded by $1\leftrightarrow 2$). In the
following we define dynamical systems, called \emph{ternary codings of rotations}, which
capture this phenomenon.

\begin{figure}
\centering
\begin{subfigure}{0.45\textwidth}
\centering
\begin{tikzpicture}[scale=.8]
    \filldraw[gray,opacity=.8] (0,0) -- (30:1.5) arc (30:90:1.5);
    \filldraw[gray,opacity=.5] (0,0) -- (160:1.5) arc (160:220:1.5);
    \draw[thick] (0,0) circle(1.5);

    \filldraw (90:1.5) circle(1pt);
    \node at (90:1.7) {\footnotesize{$0$}};
    \filldraw (30:1.5) circle(1pt);
    \node at (30:2) {\footnotesize{$\{-\alpha\}$}};
    \filldraw (160:1.5) circle(1pt);
    \node at (160:2) {\footnotesize{$\zeta-\alpha$}};
    \filldraw (220:1.5) circle(1pt);
    \node at (220:1.8) {\footnotesize{$\zeta$}};

    \node at (60:.9) {\footnotesize{$J_1$}};
    \node at (190:.9) {\footnotesize{$J_2$}};
    \node at (120:1) {\footnotesize{$J_0$}};
    \node at (270:1) {\footnotesize{$J_0$}};

    \filldraw (310:1.5) circle(1pt);
    \node at (310:1.8) {\footnotesize{$x$}};

    \filldraw (370:1.5) circle(1pt);
    \node at (370:2.2) {\footnotesize{$R_{\alpha}(x)$}};

    \draw[|->] (310:1.40) arc (310:370:1.40);
\end{tikzpicture}
\caption{}
\label{fig:RotationCodingNp2General}
\end{subfigure}\quad\quad
\begin{subfigure}{0.45\textwidth}
\centering
\begin{tikzpicture}[scale=.8]
    \filldraw[gray,opacity=.8] (0,0) -- (30:1.5) arc (30:90:1.5);
    \filldraw[gray,opacity=.5] (0,0) -- (90:1.5) arc (90:150:1.5);
    \draw[thick] (0,0) circle(1.5);

    \filldraw (90:1.5) circle(1pt);
    \node at (90:1.7) {\footnotesize{$0$}};
    \filldraw (30:1.5) circle(1pt);
    \node at (30:2) {\footnotesize{$\{-\alpha\}$}};
    \filldraw (150:1.5) circle(1pt);
    \node at (150:1.8) {\footnotesize{$\alpha$}};

    \node at (60:.9) {\footnotesize{$J_1$}};
    \node at (120:.9) {\footnotesize{$J_2$}};

    \node at (270:1) {\footnotesize{$J_0$}};

    \filldraw (310:1.5) circle(1pt);
    \node at (310:1.8) {\footnotesize{$x$}};

    \filldraw (370:1.5) circle(1pt);
    \node at (370:2.2) {\footnotesize{$R_{\alpha}(x)$}};

    \draw[|->] (310:1.40) arc (310:370:1.40);
\end{tikzpicture}
\caption{}
\label{fig:RotationCodingNp2Tight}
\end{subfigure}
\caption{An illustration of a system of codings of rotations with more than 2 intervals. In \ref{fig:RotationCodingNp2General}
we have four intervals and in \ref{fig:RotationCodingNp2Tight} three intervals. The word
$\infw u$ in \autoref{thm:characterizationNp2Complexity}, \autoref{item:item3LowComplexity} is
a coding of the orbit of some point in the system \ref{fig:RotationCodingNp2Tight}.}
\end{figure}
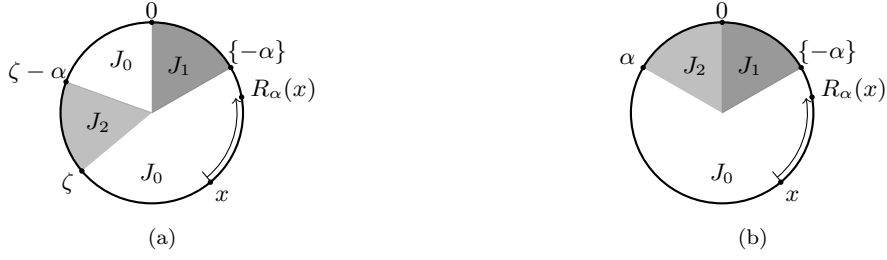

Recall the definition of Sturmian words as codings of rotations on the torus $\mathbb T$
with the half-open intervals $I_0 = I(0,1-\alpha)$ and
$I_1 = I(1-\alpha,1)$. We assume here that $\alpha< 1/2$. Take $\zeta\in I(\alpha,1-\alpha)$
and split torus $\mathbb T$ into four (three if $\zeta = \alpha$ or $\zeta = 1-\alpha$)
intervals defined by the points $0,\zeta-\alpha,\zeta$, and $1$ in increasing order:
Define the disjoint intervals $J_2 = \underline{I}(\zeta-\alpha,\zeta)$ (resp.,
$J_2 = \bar{I}(\zeta-\alpha,\zeta)$), $J_1 = I_1$, and $J_0 = I_0\setminus J_2$.
We must be careful with the value $\zeta = \alpha$ (resp., $\zeta = 1-\alpha$): If
$1\in I_1$ (resp., $1-\alpha\in I_1$) then $J_2 = \bar{I}(0,\alpha)$ (resp.,
$J_2 = \underline{I}(1-2\alpha,1-\alpha)$). Take the rotation $R_\alpha$ and
the encoding $\nu:\mathbb T\to \{0,1,2\}$, $x\mapsto i$ if and only if $x\in I_i$. The
word $\infw t_{\alpha,\zeta,x} = (\nu(R_{\alpha}^n(x))$ is called the \emph{rotation word}
of slope $\alpha$, \emph{offset} $\zeta$, and intercept $x$. See
\autoref{fig:RotationCodingNp2General} for an illustration. When indicate the choices
of endpoints of $J_i$ as follows.
If $1\in J_1$ and $\zeta \in J_2$ (resp., $1\notin J_1$, $\zeta \notin J_2$)
we denote the obtained word by $\infw t_{\overline{\alpha},\overline{\zeta},\rho}$ (resp.,
$\infw t_{\underline{\alpha},\underline{\zeta},\rho}$).
If $1\in J_1$ and $\zeta \notin J_2$ (resp., $1\notin J_1$, $\zeta \in J_2$),
we denote this by $\infw t_{\overline{\alpha},\underline{\zeta},\rho}$ (resp.,
$\infw t_{\underline{\alpha},\overline{\zeta},\rho}$).
Notice that $\infw t_{\overline\alpha,\underline{\alpha},\rho}$ and
$\infw t_{\underline{\alpha},\overline{1-\alpha},\rho}$ are not defined: this would imply that the intervals $J_1$ and $J_2$ overlap.

Observe now that, by the discussion following \autoref{ex:nonRecurrentInAbShift},
for $\infw u$ of factor complexity $n+2$ as in \autoref{thm:characterizationNp2Complexity}\eqref{item:item3LowComplexity}, we have $\infw u = \infw t_{\alpha,\alpha,\rho}$ for some $\rho \in \mathbb T$
(see \autoref{fig:RotationCodingNp2Tight}). Further, any word $\infw x\in \abclsr{\infw u}$
is of form $\infw t_{\alpha,\zeta,\rho'}$ for some $\zeta\in I(\alpha,1-\alpha)$, $\rho' \in \mathbb T$.
Our main goal is to show that $\infw t_{\alpha,\zeta,\rho}\in \abclsr{\infw u}$ for all possible
$\zeta\in [\alpha,1-\alpha]$.

Recall that Sturmian words are balanced, so that for each $n\in\N$ and for each $i=1,2$,
the set $\{|v|_i\colon v\in \mathcal L(\infw u)\}$ comprises two values (depending on $n$
and $\alpha$). We say that a factor $v$ is \emph{$1$-heavy} (resp., \emph{$2$-heavy})
if $|v|_1$ (resp., $|v|_2$) attains the larger of the two possible values. Otherwise
we say that $v$ is \emph{$1$-light} (resp., \emph{$2$-light}). If $v$ is $1$-heavy and
$2$-heavy, we say that $v$ is \emph{$1$-$2$-heavy}. Similarly, $v$ is called \emph{$1$-$2$-light}
if $v$ is $1$-light and $2$-light. We make use of the following result appearing in
\cite[part of Thm.~19]{RigoSalimovVandomme13}.
\begin{proposition}\label{prop:heavyPointOnCircle}
Let $\infw s$ be a Sturmian word of slope $\alpha$ and intercept $\rho$ and let $m\geq 1$.
Then the prefix of length $m$ of $\infw s$ is heavy if and only if
$\rho \in I(R_{\alpha}^{-m}(0),1)$. Here $I(R_{\alpha}^{-m}(0),1)$ contains the point $R_{\alpha}^{-m}(0)$
if and only if $1\notin I_1$.
\end{proposition}

We may apply \autoref{prop:heavyPointOnCircle} to determine whether a point starts with
a heavy factor of length $m$ or not. Indeed, for the letter $1$, the proposition stands as is:
$\infw t_{\alpha,\zeta,\rho}$ begins with a $1$-heavy factor of length $m$ if and only if
$\rho\in I(\{-m\alpha\},1)$. For $2$-heavy factors we take into account the
rotation induced by the offset $\zeta$. Thus the word $\infw t_{\alpha,\zeta,\rho}$
begins with a $2$-heavy factor if and only if
$\rho \in I(\{-m\alpha + \zeta\},\zeta)$.
Here the interval $I(\{-m\alpha + \zeta\},\zeta)$ contains the point $\{-m\alpha+\zeta\}$ (resp., $\zeta$)
if and only if $\zeta \notin J_2$ (resp., $\zeta \in J_2$).
We define the following distance on the torus: $\lVert x \rVert = \min\{x,1-x\}$.
Thus, e.g., $\max\{x,1-x\} = 1-\lVert x \rVert $.

\begin{lemma}\label{lem:1-2HeavyFactors}
Let $\infw x = \infw t_{\alpha,\zeta,\rho}$. Then
\begin{enumerate}[topsep=0pt]
\item $\infw x$ contains a $1$-heavy--$2$-light and a $2$-heavy--$1$-light factor for each length.

\item There exists a $1$-$2$-heavy factor
$v\in \mathcal L(\infw x)$ of length $m$ if and only if $\{-m\alpha\}<1-\lVert \zeta \rVert$,
or $\{-m\alpha\} = 1-\lVert \zeta \rVert$ and $\infw{x} = \infw{t}_{\overline{\alpha}, \underline{\{m\alpha\}},\{-n\alpha\}}$
or $\infw{x} = \infw{t}_{\underline{\alpha}, \overline{\{-m\alpha\}},\{-(m+n)\alpha\}}$ for some $n\geq 0$.

\item There exists a $1$-$2$-light factor $v\in \mathcal L(\infw x)$
of length $m$ if and only if $\{-m\alpha\}>\lVert \zeta \rVert $,
or $\{-m\alpha\} = \lVert \zeta \rVert$ and
$\infw{x} = \infw{t}_{\overline{\alpha},\underline{\{-m\alpha\}},\{-(m+n)\alpha\}}$ or $\infw{x} = \infw{t}_{\underline{\alpha},\overline{\{m\alpha\}},\{-n\alpha\}}$. 
\end{enumerate}
\end{lemma}

\begin{proof}
We give a proof case by case. Consider factors of length $m$ and write $\mu = \{-m\alpha\}$ for short.
\begin{enumerate}[leftmargin=*]
\item  We first consider $1$-heavy--$2$-light factors.
By the preceding observations on $1$-heavy and $2$-heavy
factors, $\infw{x}$ has a $1$-heavy--$2$-light factor if and only if
$I(\mu,1)\cap I(\zeta,\{\zeta+\mu\})\neq \emptyset$. The interval
$I(\max\{\zeta,\mu\},\min\{\zeta+\mu,1\})$ is always in the intersection, since
$\zeta,\mu< 1$ and $\zeta,\mu< \zeta+\mu$. Since $(\{-n\alpha\})_n$ is dense in $[0,1)$,
some shift of $\infw x$ corresponds to a coding of a point in this interval.

Similarly, $\infw{x}$ has a $2$-heavy--$1$-light factor if and only if
$I(0,\mu)\cap I(\{\zeta+\mu\},\zeta)\neq \emptyset$. The interval
$I(\max\{0,\mu+\zeta-1\},\min\{\zeta,\mu\})$ is always in the intersection, since
$\zeta, \mu > 0$ and $\zeta,\mu > \mu + \zeta-1$.

\item
We then consider $1$-$2$-heavy factors. Similar to above, $\infw{x}$ has a $1$-$2$-heavy
factor if and only if $I(\mu,1)\cap I(\{\mu+\zeta\},\zeta) \neq \emptyset$.
Assume first that $\mu<1-\lVert \zeta \rVert$.
If $\mu < \zeta$ then $I(\mu,\zeta)$ is in the intersection.
If $\zeta \leq \mu < 1-\lVert \zeta \rVert$, then $\mu+\zeta <1$
so that $\{\mu+\zeta\} = \mu+\zeta$ and $I(\mu+\zeta,1)$ is in the intersection.
The denseness of $(\{-n\alpha\})_{n\in\N}$ in $\mathbb{T}$ again implies that some shift of $\infw x$ corresponds to a point in this interval.

Assume then that $\mu = 1-\lVert \zeta \rVert$. If $\zeta=\lVert \zeta \rVert$,
then $\{\mu+\zeta\} = 0$. Now $I(\{\mu+\zeta\},\zeta)$ and $I(\mu,1)$ can
share at most one point in common, namely the point $1$. Now the intersection
is non-empty if and only if $1\in J_1$ and $\zeta\notin J_2$.
Further, $\infw t_{\overline{\alpha},\underline{\zeta},0}$ is the only word starting with a
$1$-$2$-heavy factor. To hit the point $0$ in the orbit starting from $\rho \in \mathbb{T}$, we must have $\rho = \{-n\alpha\}$ for some $n\geq 0$. So, we have $\infw{x}$ contains a $1$-$2$-heavy factor of length $m$ if and only if $\infw{x} = \infw{t}_{\overline{\alpha},\underline{\zeta},\{-n\alpha\}}$, where
$\zeta = 1 - \mu = 1 - \{-m\alpha\} = \{m\alpha\}$.

Similarly, if $\zeta = 1-\lVert \zeta \rVert$, then $I(\{\mu+\zeta\},\zeta)$ and
$I(\mu,1)$ can share at most one point in common, namely $\zeta$.
The intersection is not empty if and only if
$1\notin J_1$ and $\zeta\in J_2$. In this case $\infw t_{\underline{\alpha},\overline{\zeta},\zeta}$ is
the only point starting with a $1$-$2$-heavy factor. To hit the point $\zeta$ in the orbit of $\rho$, we must have $\rho = \{\zeta-n\alpha\}$
for some $n\geq 1$. Hence $\infw{x}$ contains a $1$-$2$-heavy factor if
and only if $\infw{x} = \infw{t}_{\underline{\alpha},\overline{\zeta},\{\zeta - n\alpha\}}$,
where $\zeta = 1 - \|\zeta\| = \mu = \{-m\alpha\}$.

Assume finally that $\mu>1-\lVert \zeta \rVert$. It follows that
$\mu > \zeta$ and $\mu+\zeta>1$ and thus $0<\{\mu+\zeta\}<\zeta$. Therefore
$I(\mu,1)\cap I(\{\mu+\zeta\},\zeta) = \emptyset$. This concludes the case of $1$-$2$-heavy
factors.

\item
Let us then finally consider $1$-$2$-light factors. We proceed analogous to the previous
case. The word $\infw{x}$ has a $1$-$2$-light factor of length $m$ exists if and only if
$I(0,\mu)\cap I(\zeta,\{\zeta+\mu\}) \neq \emptyset$. Assume first
that $\mu>\lVert \zeta \rVert $. If $\mu > \zeta$, then the interval $I(\zeta,\mu)$ is in
the intersection. If $\zeta > \mu \geq 1-\lVert \zeta \rVert $, then $\mu+\zeta > 1$ so that
$\{\mu + \zeta\}>0$. Now $I(0,\{\mu+\zeta\})$ is in the intersection.

Assume then that $\mu = \lVert \zeta \rVert$. If $\zeta=\lVert \zeta \rVert$,
then $\mu+\zeta <1$. Now $I(\zeta,\mu+\zeta)$ and $I(0,\mu)$ can
share at most one point in common, namely the point $\zeta$. By the observations preceding the lemma,
the intersection is not empty if and only if $1\in J_1$ and $\zeta\notin J_2$.
Now $\infw t_{\overline{\alpha},\underline{\zeta},\zeta}$ is the only point starting
with a $1$-$2$-light factor. The only way to hit the point $\zeta$ in the orbit of $\rho$ is that $\rho = \{\zeta - n\alpha\}$ for some $n\geq 0$. It follows that $\infw{x}$ contains a $1$-$2$-light factor if and only
if $\infw{x} = \infw{t}_{\overline{\alpha},\underline{\zeta},\{\zeta-n\alpha\}}$, where $\zeta = \{-m\alpha\}$.

Similarly, if $\zeta = 1-\lVert \zeta \rVert$, then $\mu+\zeta = 1$. Now
$I(\zeta,\{\mu+\zeta\})$ and $I(0,\mu)$ can share at most one point in common, namely the point $1$.
By the observations preceding the lemma, the intersection is not empty if and only if
$1\notin J_1$ and $\zeta\in J_2$. Now $\infw t_{\underline{\alpha},\overline{\zeta},0}$ is the only
factor starting with a $1$-$2$-light factor. Again, we have $\infw{x}$ contains a $1$-$2$-light factor if and only if $\infw{x} = \infw{t}_{\underline{\alpha},\overline{\zeta},\{-n\alpha\}}$,
where $\zeta = 1  - \mu = \{m\alpha\}$.

Finally, if $\mu<\lVert \zeta \rVert $, then $\mu < \zeta$ and $\mu+\zeta \leq \mu + 1-\lVert \zeta \rVert <1$.
Thus $I(0,\mu)$ and $I(\zeta,\mu+\zeta)$ do not intersect.
This concludes the proof.\qedhere
\end{enumerate}
\end{proof}

As is evident from \autoref{lem:1-2HeavyFactors}(2), the existence of a
$1$-$2$-heavy factor of a certain length depends not only on $\zeta$, but also on $\rho$ and how the
endpoints of the intervals are defined.
For example, the word $\infw{t}_{\overline{\alpha}, \underline{\{m\alpha\}},0}$ begins
with a $1$-$2$-heavy factor of length $m$, while
$\infw{t}_{\underline{\alpha}, \overline{\{m\alpha\}},\{-n\alpha\}}$ does
not contain such a factor. Note further that
$\infw{t}_{\overline{\alpha}, \underline{\{m\alpha\}},0}$ contains only
one occurrence of such a factor, and hence is non-recurrent.
In fact, any word
$\infw{t}_{\overline{\alpha}, \underline{\{m\alpha\}},\{-n\alpha\}}$,
$n\geq 0$, contains exactly one such factor of length $m$, namely at position $n$. 

\begin{lemma}\label{lem:largerOffsetAbelianClsr}
If $\lVert \zeta\lVert> \lVert \zeta' \rVert$ then $\infw t_{\alpha,\zeta,\rho}\in \abclsr{\infw t_{\alpha,\zeta',\rho'}}$
but $\infw t_{\alpha,\zeta',\rho'}\notin \abclsr{\infw t_{\alpha,\zeta,\rho}}$. 
\end{lemma}
\begin{proof}
Let $\infw{x} = \infw t_{\alpha,\zeta,\rho}$ and $\infw{u} = \infw t_{\alpha,\zeta',\rho'}$ for short.
Observe that for any $w,w'$, where $w\in \mathcal L_m(\infw u)$ and $w'\in \mathcal L_m(\infw x)$,
we have $||w|_1-|w'|_1|\leq 1$ and $||w|_2-|w'|_2|\leq 1$.

Let us first show that $\infw{x} \in \abclsr{\infw{u}}$.
By \autoref{lem:1-2HeavyFactors}(1), both words contain both $1$-heavy-$2$-light and $2$-heavy-$1$-light factors of each length.
We show that whenever $\infw{x}$ contains a $1$-$2$-heavy factor or a $1$-$2$-light factor length $m$, then $\infw{u}$ contains such a factor as well, which suffices for the claim. To this end, let $w\in \mathcal L_n(\infw x)$.
If $w$ is a $1$-$2$-heavy factor, then by \autoref{lem:1-2HeavyFactors}(2),
we have
$\{-m\alpha\}\leq 1-\lVert \zeta \rVert < 1-\lVert \zeta' \rVert $ so that $\infw u$ contains a $1$-$2$-heavy
factor by the same lemma. If $w$ is a $1$-$2$-light factor, then by \autoref{lem:1-2HeavyFactors}(3),
$\{-m\alpha\}\geq \lVert \zeta \rVert > \lVert \zeta' \rVert $ so that $\infw u$ again contains a $1$-$2$-light factor of length $m$.

We then show that $\infw{u} \notin \abclsr{\infw{x}}$.
Since $(\{-m\alpha\})_{m\geq 1}$ is dense in $[0,1)$, there must exist $m\in \N$ for which
$\lVert \zeta \rVert > \{-m\alpha\} > \lVert \zeta' \rVert $. By
\autoref{lem:1-2HeavyFactors}(3), $\infw u$ contains a $1$-$2$-light factor of length $m$, while $\infw x$ does not. It
follows that $\infw u\notin \abclsr{\infw x}$.
\end{proof}

We may now characterize the abelian closure of $n+2$ factor complexity words via ternary 
codings of rotations.

\begin{proposition}\label{prop:characterizationNp2Uncountable}
Let $\infw u = \infw t_{\alpha,\alpha,\rho}$ for some $\rho\in \mathbb T$. Then
$\abclsr{\infw u} =
                    \bigcup_{\zeta\in [\alpha ,1-\alpha]}
                            \soc{\infw t_{\alpha,\zeta,\rho}}.$
\end{proposition}
\begin{proof}
By the above lemma we have $\infw t_{\alpha,\zeta,\rho} \in \abclsr{\infw u}$ for all $\zeta\in(\alpha,1-\alpha)$.
For $\zeta=\alpha$ or $\zeta = 1-\alpha$, all words $\infw t_{\alpha,\zeta,\rho}$ either contain or do not contain
a $1$-$2$-light (resp. $1$-$2$-heavy) factor regardless of $\rho$. (Recall that the words
$\infw t_{\overline{\alpha},\underline{\alpha},\rho}$ and
$\infw t_{\underline{\alpha},\overline{1-\alpha},\rho}$ are not defined.)
As there are no other words in
$\abclsr{\infw{u}}$, this concludes the proof.
\end{proof}

In fact, utilising \autoref{lem:largerOffsetAbelianClsr} we can characterize
the abelian closure any word $\infw{t}_{\alpha,\zeta,\rho}$.
The proof above applied to the setting $\|\zeta\| = \|\alpha\|$, i.e., when $\infw{u}$
is a minimal complexity word, carries over to arbitrary $\zeta$ with minor
modifications:

\begin{proposition}
Let $\infw{u} = \infw{t}_{\alpha,\zeta,\rho}$ with $\|\zeta\| > \|\alpha\|$.
Then
\begin{equation*}
    \abclsr{\infw{u}} = \bigcup_{\substack{ \|\xi\| \geq \|\zeta\| \\ \rho' \in \mathbb{T}}} \soc{\infw{t}_{\alpha,\xi,\rho'}}  \setminus S,
\end{equation*}
where $S$ is a countable
set of words depending on $\zeta$ and $\rho$ as follows.

\begin{enumerate}
\item If $1- \|\zeta\|, \|\zeta\| \notin \{\{-m\alpha\} \colon m \in \N \}$,
then $S = \emptyset$.

\item Assume that $1 - \|\zeta\|  =  \{-m \alpha\}$ for some $m \geq 1$.
If $\infw{u} = \infw{t}_{\overline{\alpha}, \underline{\{m\alpha\}},\{-n\alpha\}}$
or $\infw{u} = \infw{t}_{\underline{\alpha}, \overline{\{-m\alpha\}},\{-(m+n)\alpha\}}$ for some $n\geq 0$, then $S = \emptyset$.
Otherwise
\begin{equation*}
S = \{ \infw{t}_{\overline{\alpha}, \underline{\{m\alpha\}},\{-n\alpha\}} \colon n \geq 0\} \cup \{ \infw{t}_{\underline{\alpha}, \overline{\{-m\alpha\}},\{-(m+n)\alpha\}} \colon n \geq 0\}.
\end{equation*}

\item Assume that $\|\zeta\| = \{-m\alpha\}$ for some $m\geq 1$.
If $\infw{u} = \infw{t}_{\overline{\alpha},\underline{\{-m\alpha\}},\{-(m+n)\alpha\}}$ or
$\infw{u} = \infw{t}_{\underline{\alpha},\overline{\{m\alpha\}},\{-n\alpha\}}$ for some $n\geq 0$, then $S = \emptyset$. Otherwise
\begin{equation*}
S = \{ \infw{t}_{\overline{\alpha},\underline{\{-m\alpha\}},\{-(m+n)\alpha\}} \colon n \geq 0\} \cup \{ \infw{t}_{\underline{\alpha},\overline{\{m\alpha\}},\{-n\alpha\}} \colon n \geq 0\}.
\end{equation*}
\end{enumerate}
\end{proposition}

\begin{proof}
Notice that any word $\infw{y} \in \abclsr{\infw{u}}$, we have
that $\infw{y}$ is of the form $\infw{t}_{\alpha,\xi,\rho'}$. Indeed,
using the mappings $\varphi_1$ and $\varphi_2$ as in the discussion following \autoref{ex:nonRecurrentInAbShift},
$\varphi_1(\infw{y})$ and $\varphi_2(\infw{y})$ are Sturmian
words with slope $\alpha$. We deduce that they are interleavings
of Sturmian words giving rise to the claimed form of $\infw{y}$.

\autoref{lem:largerOffsetAbelianClsr} then gives that
$\abclsr{\infw{u}}$ is a subset of
$\bigcup_{\|\xi\| \geq \|\zeta \|}\soc{\infw{t}_{\alpha,\xi,\rho'}}$, but it is possibly a proper subset. The same lemma shows that $\abclsr{\infw{u}}$ is a superset of $\bigcup_{\|\xi\| > \|\zeta \|}\soc{\infw{t}_{\alpha,\xi,\rho}}$ in any case.

Therefore, we may focus on words $\infw{y} = \infw{t}_{\alpha,\xi,\rho'}$ with offset $\xi$ having
$\|\xi\| = \|\zeta\|$. Notice that the three points are disjoint. Indeed, in 2., we assume that
$1 - \|\zeta\| = \{-m\alpha\}$, which gives $\|\zeta\| = \{m\alpha\}$.
Hence $\|\zeta\| \neq \{-m'\alpha\}$ for any $m' \geq 1$, as otherwise
$\{(m + m')\alpha\} = 0$ which would leave $\alpha$ rational.

To identify the set of words $S$ not in 
the abelian closure of $\infw{u}$, we employ \autoref{lem:1-2HeavyFactors}.

\begin{enumerate}
\item Assume that $1 - \|\zeta\|$,
$\|\zeta \| \notin \{ \{-m\alpha\} \colon m \in \N\}$. \autoref{lem:1-2HeavyFactors}(1) then states that the
existence of a $1$-$2$-heavy or light factor does not depend on
the point whose orbit we encode, nor the choices of the endpoints of the intervals. That is to say, all words with offset $\xi$,
$\|\xi\| = \|\zeta\|$, simultaneously either have or do not
have a $1$-$2$-heavy (resp., light) factor of length $m$ independent to the
choice of starting point $\rho'$ of the orbit. This suffices to show that
$S = \emptyset$ in this case.

\item
Assume that $1 - \|\zeta\| = \{-m\alpha\}$ for some $m\geq 1$. There is only
one length of factors in which the existence of a $1$-$2$-heavy factor
depends on the starting point $\rho$ and the choice of the endpoints of the
intervals. This length is $m$. By \autoref{lem:1-2HeavyFactors}(2),
if $\infw{u} = \infw{t}_{\overline{\alpha}, \underline{\{m\alpha\}},\{-n\alpha\}}$
or $\infw{u} = \infw{t}_{\underline{\alpha}, \overline{\{-m\alpha\}},\{-(m+n)\alpha\}}$ for some $n\geq 0$, then the word
contains such a factor. In this case $S = \emptyset$. If $\infw{u}$ is not of this form, then it does not contain such a factor, while all the words
$\infw{t}_{\overline{\alpha}, \underline{\{m\alpha\}},\{-n\alpha\}}$ and
$\infw{t}_{\underline{\alpha}, \overline{\{-m\alpha\}},\{-(m+n)\alpha\}}$, $n\geq 0$, do. The claim then follows.

\item This is analogous to the one above.\qedhere
\end{enumerate}
\end{proof}

\subsection{Recurrent minimal complexity words with at least four letters}

Surprisingly, for alphabet of size greater than 3 there are always
only finitely many subshifts:

\begin{theorem} \label{th:min_compl}
Let $ u$ be a recurrent word of factor complexity $n+C$ for all
$n\geq 1$, where $C>2$. Then $\abclsr{\infw u}$ contains exactly two
minimal subshifts.
\end{theorem}

The proof is based on the characterization of words of factor
complexity $n+C$ for all $n\geq 1$ from \cite{FM97}.

\begin{lemma}[{\cite[Lem.~4]{FM97}}]\label{lem:min_compl} Let $\infw u$ be a recurrent word of
minimal complexity an alphabet $A$; then there exist distinct elements $e_1, \dots, e_b$,
$f_1, \dots, f_c$, $g_1, \dots, g_d$ in $A$ such that the sets $E=\{e_1, \dots, e_b\}$,
$F=\{f_1, \dots, f_c\}$, and $ G=\{g_1, \dots, g_d\}$ are pairwise disjoint,
$E \cup F \cup G = A$, with $G \neq \emptyset$, and $E \cup F \neq \emptyset$, and there
exists a Sturmian word $\infw{s}$ on $\{0, 1\}$ such that, if $\sigma$ is the
substitution
\begin{equation*}
    \begin{cases}0 &\mapsto g_1 \cdots g_d e_1 \cdots e_b\\
    1 &\mapsto g_1 \cdots g_d f_1 \cdots f_c
    \end{cases},
\end{equation*} then $\sigma(\infw s)=W \infw u$, where $W$ is a (possibly
empty) prefix of $\sigma(0)$ or $\sigma(1)$.
\end{lemma}

\begin{lemma}\label{lem:blocks} Let $\infw w$ be a recurrent word of of minimal complexity an alphabet $A$ of cardinality at least 3. Then each $\infw w'\in \abclsr{\infw w}$ is a concatenation of blocks of the form $\sigma(0)$ and $\sigma(1)$ (or their reversals), where $\sigma$ is as in the previous lemma, preceded by a possibly empty suffix of $\sigma(0)$ or $\sigma(1)$ (or a reversal of a prefix).   
\end{lemma}

\begin{proof}
The proof is quite direct. Since $|A|\geq 4$, at least one of the sets $G,E,F$ contains at least two letters. Let it be $E$ (for other sets it is similar). First we show that for any $\infw w'\in \abclsr{\infw w}$ the letters from $E$ must occur in blocks $e_1\cdots e_b$ (or in $e_b\cdots e_1$ -- this case is symmetric, all the blocks are reversed). For this, it is enough to consider only factors of length 2 and 3. Indeed, if $e_1$ occurs in $\infw w'$, then the only factors containing $e_1$ in $\infw w$ are $e_1 e_2$ and $g_d e_1$. With the exception of the case $F= \emptyset$ and $|G|=1$, we have that $g_de_1g_d$ is not an abelian factor of $\infw w$. Since $g_de_1g_d$ is not an abelian factor of $\infw w$, in $\infw w'$ we must have $g_d e_1 e_2$ or $e_2 e_1 g_d$. Continuing this line of reasoning with $e_2, e_1, e_3$ instead of $e_1, g_d, e_2$ etc., we get that $\infw w'$ the letters from $E$ must occur in blocks $e_1\cdots e_b$ (or in $e_b\cdots e_1$). 
If $F= \emptyset$ and $|G|=1$, then $|E|\geq 3$ (the cardinality of the alphabet is at least 4), and we can start by  $e_2$: by considering factors of length 2 and 3, we see that it can occur only in factors $e_1 e_2 e_3$ and $e_3 e_2 e_1$. The rest of the proof is the same.

In the same way we prove that each such block  $e_1\cdots e_b$ must be surrounded by $g_1 \cdots g_d$ from both sides, so we have $\sigma(0)g_1\cdots g_d$. 
Now we show that this block $\sigma(0)$ must be followed by either $\sigma(0)$ or $\sigma(1)$. We already have the beginning of the block $g_1\dots g_d$. After it, one must have either $e_1$, or $f_1$, or $g_1$ in the case if $F=\emptyset$  (again, it is enough to consider factors of length 2 and 3 containing $g_d$). In the cases of $e_1$ or $f_1$ it must continue with $e_2\cdots e_b$ or $f_2\cdots f_c$, respectively, thus finishing the block $\sigma(0)$ or $\sigma(1)$. In the case when $F$ is empty, we already have $\sigma(1)$. In the same way one can show that after the block $\sigma(1)$ one must also have a full block $\sigma(0)$ or $\sigma(1)$. 
\end{proof}

We remark that the cardinality at least 4 of the alphabet is essential for the above lemma. In the case $F= \emptyset$, $|G|=1$ and $|E|=2$ the letters $e_1$ and $e_2$ can be separated by $g_d$, which corresponds to Theorem \ref{thm:abelianSSMinimalComplexity} \eqref{item:item3LowComplexity}, when we have uncountably many minimal subshifts. The fourth letter blocks this possibility: either we have $|E|\geq3$, in which case $e_2$ ``glues'' letters $e_1$ and $e_2$, or $|G|\geq 2$, so the two letters $g_1$ and $g_2$ prevent mixing. \medskip

Let $\infw{u}$ be uniformly recurrent. We define a word
$\infw{u}^R$ for which $\lang{\infw{u}^R} = \lang{\infw{u}}^R$,
i.e., the set of reversals of the factors of $\infw{u}$.
Indeed, take the sequence $(p_n)_n$ of prefixes of $\infw{u}$,
and consider the sequence $(p_n^R)$ of their reversals.
There is a subsequence which converges to an infinite word $\infw{u}^R$. We claim that
$\lang{\infw{u}^R} = \lang{\infw{u}}^R$. As $\infw{u}^R$ is constructed using reversals of factors of $\infw{u}$, we have
that $\lang{\infw{u}^R} \subseteq \lang{\infw{u}}$. Let $x \in \lang{\infw{u}}$. Since $\infw{u}$ is uniformly recurrent,
$x$ must occur in $p_n$ for $n$ large enough. As $x$ occurs
within bounded gaps, we conclude that the words in the converging
subsequence of $(p_n)_{n}$ must have $x^R$ occurring for all $n$
large enough, the first occurrence occurring with a uniform bound. Hence $x^R \in \lang{\infw{u}^R}$.

Notice that we immediately have that
$\infw{u}^R \in \abclsr{\infw{u}}$.

\begin{proof}[Proof of \autoref{th:min_compl}.]
The two subshifts are $\soc{\infw{u}}$ and $\soc{\infw{u}^R}$. Suppose that there exists a word
$\infw{w} \in \abclsr{\infw{u}}$ such that it is not from $\soc{\infw{u}}$. Due to \autoref{lem:blocks}, cutting a
short prefix of $\infw{w}$, we get a word $\infw{w}'$ such that $\infw{w}'=\sigma(\infw{v}')$,
$\infw{v}'\in\{0,1\}^{\mathbb{N}}$, and $\infw{v}'$ is not in the
shift orbit closure of $\infw{v}$, where $\infw{v}$ is
a Sturmian word as in \autoref{lem:min_compl}. So, $\infw{v}'$ contains a factor $w'$ which is not
abelian equivalent to any factor of $\infw{v}$. It is straightforward to see that then $\sigma(w')$ is not
abelian equivalent to any factor of $\infw{u}$. Indeed, it cannot be abelian equivalent to a factor
consisting of full blocks. And if is happens to be abelian equivalent to a factor which does
not consist of full blocks, then it is also equivalent to a shift of this factor that
consists of full block, which is not possible. So, $\sigma(w')$ is not an abelian factor of
$\infw{u}$, hence $\infw{u}'$ is not in $\abclsr{\infw{u}}$, a contradiction.
\end{proof}

\section{Abelian closures of Arnoux--Rauzy words} \label{sec:AR}

In this section we discuss Arnoux--Rauzy words, which
are another generalization of Sturmian words to larger alphabet.
One of the ways to define Arnoux--Rauzy words is via palindromic
closures. The following basics on Arnoux--Rauzy words are
well-known and mostly taken from~\cite{ArRa,DJP}. In fact, this is
a generalization of the facts about Sturmian words given for
binary words in~\cite{deLuca}.

A finite word $v=v_0\cdots v_{n-1}$ is a \emph{palindrome} if it
is equal to its reversal, i.e., $v=v_{n-1}\cdots v_0$.
The \emph{right palindromic closure} of a finite word $u$, denoted
by $u^{(+)}$, is the shortest palindrome that has $u$ as a prefix.
The \emph{iterated (right) palindromic closure operator} $\psi$ is
 defined recursively by the following rules:
\[ \psi(\varepsilon)=\varepsilon, \quad \psi(va)=(\psi(v)a)^{(+)}\]
for all $v \in  \Sigma^*$ and $a \in  \Sigma$.
The definition of $\psi$ may
be extended to infinite words $\infw u$ over $\Sigma$ as
$\psi(\infw u)=\lim_{n} \psi(\pref_n (\infw u))$, i.e., $\psi(\infw u)$ is the
infinite word having $\psi(\pref_n(\infw u))$ as its prefix for
every $n \in \N$.

Let $\Delta$ be an infinite word on the alphabet $\Sigma$ such
that every letter occurs infinitely often in $\Delta$. The word
$\infw{c} = \psi(\Delta)$ is then called a \emph{characteristic (or
standard) Arnoux--Rauzy word} and $\Delta$ is called the
\emph{directive sequence} of $\infw{c}$. An infinite word $\infw{u}$ is called
an Arnoux--Rauzy word if it has the same set of factors as a
(unique) characteristic Arnoux--Rauzy word, which is called the
characteristic word of $\infw{u}$. The directive sequence of an
Arnoux--Rauzy word is the directive sequence of its characteristic
word. An example of Arnoux--Rauzy word is given by the Tribonacci word $\infw T$, which can be defined as the fixed point of the morphism $0\to 01$, $1\to 02$, $2 \to 0$. It is not hard to see that the Tribonacci word is an Arnoux--Rauzy word with the directive
sequence $(012)^\omega$.

Apparently, the structure of abelian closures of Arnoux--Rauzy
words is rather complicated. For example, it is not hard to see
that for any Arnoux--Rauzy word with a characteristic word $\infw{c}$ its
abelian closure contains $20\infw{c}$ (here we assume that $0$ is the
first letter of $\Delta$ and $2$ is the third letter occurring in
$\Delta$ for the first time, i.e., $\Delta$ has a prefix of the
form $0\{0,1\}^*1 \{0,1\}^* 2$). On the other hand,
$20\infw{c} \notin \soc{\infw{c}}$.

T. Hejda, W. Steiner, and L.Q. Zamboni studied the abelian shift
of the Tribonacci word $\infw{T}$. They announced that $\mathcal
A_{\infw T}\setminus \soc{\infw{T}}\neq\emptyset$ but that
$\soc{\infw{T}}$ is the only minimal subshift contained in
$\abclsr{\infw{T}}$
\cite{HejdaSteinerZamboni15, ZamboniPersonal}.

An interesting open question is to understand the general
structure of Arnoux--Rauzy words (see Problem \ref{prob:AR}).

\section{Abelian closures of general subshifts} \label{sec:general_subshifts}

In this paper and the previous works on the topic, the focus has been on abelian
closures of infinite words. It would be interesting to investigate properties of abelian
closures of general subshifts.

We recall some definitions from \cite[\S1.5]{MR1905123}, but slightly modify the terminology.
A bi-infinite word $\infw{x} \in \Sigma^{\Z}$ is said to \emph{avoid} a set of words
$\mathcal{F} \subseteq \words$ if $\mathcal{L}(\infw{x}) \cap \mathcal{F} = \emptyset$. Let
$X_{\mathcal{F}}$ denote the set of bi-infinite words avoiding $\mathcal{F}$. A
\emph{subshift} is a set $X_{\mathcal{F}}$ for some $\mathcal{F}$. The shift operator is
defined similar to the case of infinite words. Now a set $X \subseteq \Sigma^{\Z}$ is a
subshift if and only if $\sigma(X) = X$ and is closed in the usual topology on bi-infinite
words.

Let $X$ be a subshift and let $I(X) = \words \setminus \mathcal{L}(X)$. Define the set
$\mathcal{F}(X)$ as the set of elements of $I(X)$ which are minimal for the factor ordering,
i.e., have no proper factor in $I(X)$. Then $X = X_{\mathcal{F}(X)}$. 

\begin{definition}
If a subshift $X = X_{\mathcal{F}}$ for some finite set $\mathcal{F} \subseteq \words$,
then $X$ is called a \emph{subshift of finite type} (SFT).
If, on the other hand, $\mathcal{F}$ can be taken regular, then $X$ is called \emph{sofic}.
\end{definition}
A set $X$ is a SFT if and only if $\mathcal{F}(X)$ is finite. Similarly, $X$ is sofic if and
only if $\mathcal{F}(X)$ is regular.

We may define the abelian closure of a subshift straightforwardly.

\begin{definition}
Let $X$ be a subshift. Then its \emph{abelian closure} $\abclsr{X}$ is defined as
$\cup_{\infw{x} \in X}\abclsr{\infw{x}}$.
\end{definition}

We remark that in the previous text we considered one-way infinite words, as more customary in combinatorics on words, whereas here for general subshifts it is more natural to consider bi-infinite words. Actually, there is no principal difference for our considerations, as all the results can easily be reformulated for one-way or two-way infinite words.

We conclude this paper with a couple of examples of abelian closures of subshifts.
\begin{example}
Clearly $\abclsr{\Sigma^{\Z}} = \Sigma^{\Z}$. Let $\mathcal{F} = \{11\} \subseteq \{0,1\}^*$ and set $X = X_{\mathcal{F}}$. The subshift $X \subseteq \{0,1\}^{\Z}$ is called
the \emph{golden mean} subshift. Consider the abelian closure of $X_{\mathcal{F}}$: it
comprises those words for which all $1$s are isolated. But this is just $X_{\mathcal{F}}$
itself. Thus $\abclsr{ X_{\mathcal{F}} } = X_{\mathcal{F}}$.
\end{example}

In the above example, both subshifts are of finite type. It was concluded that they are, in
fact, their own abelian closures. This property is of course not general for SFTs, as is
shown by the following example. In fact, the abelian closure of a SFT is not in general a
SFT. 

\begin{example}
Consider the SFT $X = X_{\mathcal{F}}$ with $\mathcal{F} = \{aa,ac,ba,bb,cb\}$.
It can be characterized as the set of two-way infinite walks on the following graph.
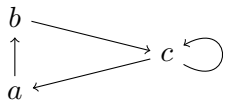
\begin{center}
\begin{tikzpicture}
\node (a) at (0,0) {$a$};
\node (b) at (0,1) {$b$};
\node (c) at (2,.5) {$c$};

\draw (a) edge[->] (b);
\draw (b) edge[->] (c);
\draw (c) edge[->] (a);
\draw (c) edge[->,out=330,in=30,looseness=8] (c);
\end{tikzpicture}
\end{center}
Assume for a contradiction, that $\abclsr{X}$ is a SFT, with
$\mathcal{F}(X) = \mathcal{F}'$. There is an integer $n$ for which each element of
$\mathcal{F}'$ has length at most $n$. Consider the word
$\infw{x} = {}^{\omega}\!c \cdot ab \cdot c^n \cdot ba \cdot c^{\omega}$. Here for a finite word $v$ by ${}^{\omega}v$ we mean the left-infinite word obtained by repeating $v$ infinitely many times. Observe that the
factors of length at most $n$ of this word occur either in
$^{\omega}\!c \cdot ba \cdot c^{\omega}$ or in $^{\omega}\!c \cdot ab \cdot c^{\omega}$. Both
of these words are in $\abclsr{X}$ by inspection, so none of the factors can be in
$\mathcal{F}'$. Thus $\infw{x}$ avoids all the forbidden factors. But,
$\infw{x} \notin \abclsr{X}$, as it contains the factor $bc^nb$. Indeed, any word in the
language $\mathcal{L}(\abclsr{X})$ that contains two occurrences of $b$ must contain at least
one occurrence of $a$.
\end{example}

The next example shows that this is also possible for binary alphabet:

\begin{example} Consider an SFT giving words of the form 

$$ \cdots 001100110011000111000111000111 \cdots,$$ 
plus $(0011)^\omega$ and $(000111)^\omega)$. It is a SFT of order 6, and its Rauzy graph contains two cycles with the same frequencies of letters and a one-way path between them.

It is readily verified that words of the form
$$ \cdots 0011001100110001100110011\cdots$$  are in the abelian closure, whereas  words of the form $$ \cdots 0011001100011(0011)^+0001100110011\cdots$$  are not. Similarly to the previous example this implies that the abelian closure is not a SFT.
\end{example}

Next we show that the abelian closure of a sofic shift is not necessarily sofic.
\begin{example}
Let $\alphabet = \{a,b,c,d\}$ be the underlying alphabet. Set 
\begin{equation*}
\mathcal{F} = \{a,b,d\} c \cup d\{a,b,c\} \cup c R d,
\end{equation*}
where $R = \{a,b\}^* \setminus (ab)^*$ and let $X = X_{\mathcal{F}}$.
Hence $X$ is of the form
\begin{equation*}
X = \{a,b\}^{\Z} \cup \{{}^{\omega}\!c \infw{x}, \infw{x}^R d^{\omega} \colon \infw{x}\in \{a,b\}^{\N}\} \cup \{{}^{\omega}\!c (ab)^n d^{\omega} \colon n\geq 0 \} \cup\{{}^{\omega}\!c^{\omega}, {}^{\omega}\!d^{\omega}\}.
\end{equation*}
(Here $\infw{x}^R$ is the left-infinite word defined by $\infw{x}$, i.e., the letter at position $-n$ of $\infw{x}^R$
equals the $n$th letter of $\infw{x}$.)

Let $\mathcal{F}' = \mathcal{F}(\abclsr{X})$. We show that
\begin{equation*}
\mathcal{F}' \cap c \{a, b\}^* d = \{cwd \colon |w|_a \neq |w|_b \}.
\end{equation*}

It follows that $\mathcal{F}'$ is cannot be regular, as the language above is well-known to
be non-regular. Hence $\abclsr{X}$ is not sofic.

Let us show the $\subseteq$ direction. Let $w$ have $|w|_a = |w|_b$. We show that the
word ${}^{\omega}\!c w d^{\omega}$ is in the abelian closure, and thus
$c w d \notin \mathcal{F}'$. Now all factors of the form $c^nx$ or $yd^m$, for $x$ a prefix
and $y$ a suffix of $w$, occur in the words ${}^{\omega}\!c w^{\omega}$ or
${}^{\omega}w d^{\omega}$, which are elements of $X$. We may thus concentrate on factors of
the form $c^n w d^m$. Now $c^n w d^m$ is abelian equivalent to $c^n (ab)^{|w|/2} d^m$ which
is clearly in the language of $X$.

Let then $w \in \{a,b\}^*$ be such that $|w|_a \neq |w|_b$. Observe that the proper factors of
$cwd$ are in the language of $X$. This means that either $cwd \in \mathcal{F}'$ or
$cwd \in \lang{\abclsr{X}}$. Assume, for a contradiction that
$cw d \in \lang{\abclsr{X}}$. Now $\Psi(cwd) = (|w|_a,|w|_b,1,1)$, but, in $\lang{X}$, any word
with Parikh vector with last two components equal to $1$ is of the form $(m,m,1,1)$. Since
$|w|_a \neq |w|_b$, there is no word in $\lang{X}$ which is abelian equivalent, so $cwd$ is not
an element of $\lang{\abclsr{X}}$. We conclude that $cwd \in \mathcal{F}'$.
\end{example}

An interesting open question is to find out whether the abelian closure of a subshift of finite type always sofic (see \autoref{problem:SFT_sofic}).

\section{Conclusions}
In this paper, we introduced and studied a notion of abelian
subshifts of infinite words. The main open problem we would like to state in this paper is the following:

\begin{problem}
Characterize words for which $\abclsr{\infw x} = \soc{\infw
x}$. \end{problem} Among binary uniformly recurrent words, this property gives a
characterization of Sturmian words, but the characterization does
not extend to usual generalizations of Sturmian words to
non-binary alphabets: neither for balanced words, nor for words of
minimal complexity, nor for Arnoux--Rauzy words. A modification of
this question is to characterize words for which $\abclsr{\infw
x}$ contains exactly one minimal subshift.

For Arnoux--Rauzy words, we showed that $\abclsr{\infw x} \neq \soc{\infw 
x}$, but their abelian closure seems to have rather complicated structure, in particular, it always contains non-recurrent words. An interesting open question is to understand the general structure of Arnoux--Rauzy words:

\begin{problem} \label{prob:AR}
Characterize abelian closures of Arnoux--Rauzy words.\end{problem}

Finally, we propose the following open question about general abelian subshifts: 

\begin{problem}\label{problem:SFT_sofic}
Is the abelian closure of an SFT always sofic?
\end{problem}

\section*{Acknowledgements}
We are grateful to Joonatan Jalonen, Ville Salo, and Luca Zamboni for fruitful
discussions and helpful comments. 

	Svetlana Puzynina is partially supported by Russian Foundation of Basic Research
(grant 20-01-00488) and by the Foundation for the Advancement of Theoretical Physics and Mathematics
BASIS". Part of the research was performed while Markus Whiteland was at the Department of Mathematics and Statistics, University of Turku, Finland.

\bibliographystyle{abbrvnat}
\biboptions{sort&compress}
\bibliography{bibliography}


\end{document}